\documentclass[journal]{IEEEtran}

\usepackage{cite}
\usepackage{amsmath,bm}
\usepackage{amsthm}
\usepackage{amssymb}
\usepackage{graphicx}
\usepackage[english]{babel}
\usepackage{algorithmic}
\usepackage[linesnumbered,ruled,vlined]{algorithm2e}
\usepackage{xcolor}
\usepackage{caption}
\usepackage{mathtools,amssymb,lipsum}
\usepackage{cuted}
\usepackage{arydshln}
\usepackage{pmat}
\usepackage{dblfloatfix}
\usepackage{xfrac}
\usepackage{float}
\usepackage{slashbox}

\usepackage{blkarray}
\usepackage{mathtools}
\usepackage{multirow}

\usepackage{tikz}
\usepackage{tkz-graph}
\usepackage{graphicx}      
\usepackage{pgfplots} 
\pgfplotsset{compat=newest} 
\pgfplotsset{plot coordinates/math parser=false}
\newlength\figureheight
\newlength\figurewidth

\usepackage[utf8]{inputenc}
\usepackage{apptools}

\DeclareMathOperator{\im}{im}

\let\leq\leqslant
\let\geq\geqslant

\newcommand{\calC}{\ensuremath{\mathcal{C}}}

\newcommand{\calH}{\ensuremath{\mathcal{H}}}

\newcommand{\calM}{\ensuremath{\mathcal{M}}}
\newcommand{\calN}{\ensuremath{\mathcal{N}}}
\newcommand{\calO}{\ensuremath{\mathcal{O}}}

\newcommand{\calQ}{\ensuremath{\mathcal{Q}}}

\newcommand{\hatp}{\ensuremath{\hat{p}}}
\newcommand{\hatq}{\ensuremath{\hat{q}}}

\newcommand{\hatx}{\ensuremath{\hat{x}}}
\newcommand{\haty}{\ensuremath{\hat{y}}}

\newcommand{\hatA}{\ensuremath{\hat{A}}}
\newcommand{\hatB}{\ensuremath{\hat{B}}}
\newcommand{\hatC}{\ensuremath{\hat{C}}}
\newcommand{\hatD}{\ensuremath{\hat{D}}}

\newcommand{\hatJ}{\ensuremath{\hat{J}}}

\newcommand{\hatV}{\ensuremath{\hat{V}}}
\newcommand{\hatW}{\ensuremath{\hat{W}}}

\newcommand{\hatZ}{\ensuremath{\hat{Z}}}

\newcommand{\barS}{\ensuremath{\bar{S}}}

\newcommand{\barZ}{\ensuremath{\bar{Z}}}

\newcommand{\bmat}{\begin{matrix}}
\newcommand{\emat}{\end{matrix}}
\newcommand{\bbm}{\begin{bmatrix}}
\newcommand{\ebm}{\end{bmatrix}}
\newcommand{\bbma}{\begin{bmatrix*}[r]}
\newcommand{\ebma}{\end{bmatrix*}}
\newcommand{\bpm}{\begin{pmatrix}}
\newcommand{\epm}{\end{pmatrix}}
\newcommand{\bvm}{\begin{vmatrix}}
\newcommand{\evm}{\end{vmatrix}}
\newcommand{\bse}{\begin{subequations}}
\newcommand{\ese}{\end{subequations}}
\newcommand{\beq}{\begin{equation}}
\newcommand{\eeq}{\end{equation}}
\newcommand{\ben}{\renewcommand{\labelenumi}{\arabic{enumi}.}
\renewcommand{\theenumi}{\arabic{enumi}}\begin{enumerate}}

\newcommand{\een}{\end{enumerate}}

\newcommand{\beni}{\renewcommand{\labelenumi}{\roman{enumi}.}
\renewcommand{\theenumi}{\roman{enumi}}\begin{enumerate}}

\newcommand{\eeni}{\end{enumerate}}

\newcommand{\bena}{\renewcommand{\labelenumi}{\alph{enumi}.}
\renewcommand{\theenumi}{\alph{enumi}}\begin{enumerate}}

\newcommand{\eena}{\end{enumerate}}

\newcommand{\bit}{\begin{itemize}}
\newcommand{\eit}{\end{itemize}}
\newcommand{\bthe}{\begin{theorem}}
\newcommand{\ethe}{\end{theorem}}
\newcommand{\blem}{\begin{lemma}}
\newcommand{\elem}{\end{lemma}}
\newcommand{\bprop}{\begin{proposition}}
\newcommand{\eprop}{\end{proposition}}
\newcommand{\bex}{\begin{example}}
\newcommand{\eex}{\end{example}}
\newcommand{\bas}{\begin{assumption}}
\newcommand{\eas}{\end{assumption}}
\newcommand{\bre}{\begin{remark}}
\newcommand{\ere}{\end{remark}}
\newcommand{\bcor}{\begin{corollary}}
\newcommand{\ecor}{\end{corollary}}
\newcommand{\bdfn}{\begin{definition}}
\newcommand{\edfn}{\end{definition}}
\newcommand{\bcon}{\begin{conjecture}}
\newcommand{\econ}{\end{conjecture}}

\newcommand{\inv}{\ensuremath{^{-1}}}

\newcommand{\set}[2]{\ensuremath{\left\{#1: #2\right\}}}

\newcommand{\R}{\ensuremath{\mathbb R}}

\newcommand{\Z}{\ensuremath{\mathbb Z}}

\newcommand{\HHH}{\mathcal{H}}

\DeclareMathOperator{\blkdiag}{blkdiag}
\DeclareMathOperator{\In}{In}

\newcommand{\systr}{\begin{bmatrix}
I & 0 \\ 0 & I \\ A^\top  & C^\top  \\ B^\top  & D^\top 
\end{bmatrix}}
\newcommand{\systrr}{\begin{bmatrix}
		I & 0 \\ 0 & I \\ A  & B  \\ C  & D 
\end{bmatrix}}

\newcommand{\systrred}{\begin{bmatrix} I & 0 \\ 0 & I \\ \hat{A}^\top & \hat{C}^\top \\ \hat{B}^\top & \hat{D}^\top \end{bmatrix}}

\newcounter{todocounter}

\usepackage[prependcaption,colorinlistoftodos]{todonotes}

\newtheorem{theorem}{Theorem}
\newtheorem{lemma}{Lemma}
\newtheorem{corollary}{Corollary}
\newtheorem{proposition}{Proposition}

\theoremstyle{definition}
\newtheorem{problem}{Problem}
\newtheorem{definition}{Definition}
\newtheorem{example}{Example}
\newtheorem{assumption}{Assumption}

\theoremstyle{remark}
\newtheorem{remark}{Remark}

\newtheoremstyle{lemappstyle}
{}{}{\itshape}{}{\bfseries}{.}{.5em}{{\thmname{#1 }}{\thmnumber{#2}}{\thmnote{ (#3)}}}
\theoremstyle{lemappstyle}

\newcounter{clemapp} \counterwithin{clemapp}{section}
\newtheorem{lemapp}[clemapp]{Lemma}
\newcommand{\blemp}{\begin{lemapp}}
\newcommand{\elemp}{\end{lemapp}}

\newcounter{cthmp} \counterwithin{cthmp}{section}
\newtheorem{theoremp}[cthmp]{Theorem}
\newcommand{\bthep}{\begin{theoremp}}
\newcommand{\ethep}{\end{theoremp}}

\makeatletter
\def\endthebibliography{%
	\def\@noitemerr{\@latex@warning{Empty `thebibliography' environment}}%
	\endlist
}

\usepackage[american]{circuitikz}
\makeatletter
\ctikzset{bipoles/vsources/minus/sign/rotation/.initial=0}
\pgfcircdeclarebipole{}{\ctikzvalof{bipoles/vsourceam/height}}{vsourceAM}{\ctikzvalof{bipoles/vsourceam/height}}{\ctikzvalof{bipoles/vsourceam/width}}{
	\pgfsetlinewidth{\pgfkeysvalueof{/tikz/circuitikz/bipoles/thickness}\pgfstartlinewidth}
	\pgfpathellipse{\pgfpointorigin}{\pgfpoint{0}{\pgf@circ@res@up}}{\pgfpoint{\pgf@circ@res@left}{0}}
	\pgfusepath{draw}
	\ifpgf@circ@oldvoltagedirection
	\pgftext[bottom,rotate=90,y=\ctikzvalof{bipoles/vsourceam/margin}\pgf@circ@res@down]{$+$}
	\pgftext[top,rotate=90,y=\ctikzvalof{bipoles/vsourceam/margin}\pgf@circ@res@up]{$-$}
	\else
	\pgftext[bottom,rotate=90,y=\ctikzvalof{bipoles/vsourceam/margin}\pgf@circ@res@down]{\rotatebox{\ctikzvalof{bipoles/vsources/minus/sign/rotation}}{$-$}}
	\pgftext[top,rotate=90,y=\ctikzvalof{bipoles/vsourceam/margin}\pgf@circ@res@up]{$+$}
	\fi
}

\pgfcircdeclarebipole{}{\ctikzvalof{bipoles/cvsourceam/height}}{cvsourceAM}{\ctikzvalof{bipoles/cvsourceam/height}}{\ctikzvalof{bipoles/cvsourceam/width}}{
	\pgfsetlinewidth{\pgfkeysvalueof{/tikz/circuitikz/bipoles/thickness}\pgfstartlinewidth}
	\pgfpathmoveto{\pgfpoint{\pgf@circ@res@left}{\pgf@circ@res@zero}}
	\pgfpathlineto{\pgfpoint{\pgf@circ@res@zero}{\pgf@circ@res@up}}
	\pgfpathlineto{\pgfpoint{\pgf@circ@res@right}{\pgf@circ@res@zero}}
	\pgfpathlineto{\pgfpoint{\pgf@circ@res@zero}{\pgf@circ@res@down}}
	\pgfpathlineto{\pgfpoint{\pgf@circ@res@left}{\pgf@circ@res@zero}}
	\pgfusepath{draw}
	
	\ifpgf@circ@oldvoltagedirection
	\pgftext[bottom,rotate=90,y=\ctikzvalof{bipoles/cvsourceam/margin}\pgf@circ@res@left]{$+$}
	\pgftext[top,rotate=90,y=\ctikzvalof{bipoles/cvsourceam/margin}\pgf@circ@res@right]{$-$}
	\else
	\pgftext[bottom,rotate=90,y=\ctikzvalof{bipoles/cvsourceam/margin}\pgf@circ@res@left]{\rotatebox{\ctikzvalof{bipoles/vsources/minus/sign/rotation}}{$-$}}
	\pgftext[top,rotate=90,y=\ctikzvalof{bipoles/cvsourceam/margin}\pgf@circ@res@right]{$+$}
	\fi
}
\makeatother

\newcommand{\sig}{\mathbf{\Sigma}}
\newcommand{\dataUXY}{(U_-,X,Y_-)}
\newcommand{\sigred}{\mathbf{\Sigma}_{\hatV,\hatW}}
\newcommand{\sigredbt}{\hat{\mathbf{\Sigma}}}

\begin{document}
%
\title{From data to reduced-order models via generalized balanced truncation}
%
%
%

\author{Azka~M.~Burohman,
	Bart~Besselink, 
	Jacquelien~M.~A.~Scherpen, 
	M.~Kanat~Camlibel, 
	
\thanks{Azka Muji Burohman, Bart Besselink and Kanat Camlibel are with the Bernoulli Institute for Mathematics, Computer Science and Artificial Intelligence, University of Groningen, Nijenborgh 9, 9747 AG, Groningen, The Netherlands. Azka Muji Burohman is also with the Engineering and Technology Institute Groningen (ENTEG), University of Groningen, Nijenborgh~4, 9747 AG, Groningen, The Netherlands.
	Jacquelien M. A. Scherpen is with the Engineering and Technology Institute Groningen (ENTEG), University of Groningen, Groningen, The Netherlands.
	The authors are also with the Jan C. Willems Center for Systems and Control, Faculty of Science and Engineering, University of Groningen, The Netherlands (e-mail: a.m.burohman@rug.nl, b.besselink@rug.nl, j.m.a.scherpen@rug.nl, m.k.camlibel@rug.nl)}
\thanks{This paper is based on research developed in the DSSC Doctoral Training Programme, co-funded through a Marie Sklodowska-Curie COFUND (DSSC 754315).}
}

\maketitle

\begin{abstract}
This paper proposes a data-driven model reduction approach on the basis of noisy data.  Firstly, the concept of data reduction is introduced. In particular, we show that the set of reduced-order models obtained by applying a Petrov-Galerkin projection to all systems explaining the data characterized in a large-dimensional quadratic matrix inequality (QMI) can again be characterized in a lower-dimensional QMI. Next, we develop a data-driven generalized balanced truncation method that relies on two steps. First, we provide necessary and sufficient conditions such that systems explaining the data have common generalized Gramians. Second, these common generalized Gramians are used to construct projection matrices that allow to characterize a class of reduced-order models via generalized balanced truncation in terms of a lower-dimensional QMI by applying the data reduction concept.
Additionally, we present alternative procedures to compute a priori and a posteriori upper bounds with respect to the true system generating the data.
Finally, the proposed techniques are illustrated by means of application to an example of a system of a cart with a double-pendulum. 
\end{abstract}


%
\IEEEpeerreviewmaketitle

\section{Introduction}
Model reduction refers to the problem of constructing low-dimensional system models that accurately approximate complex high-dimensional systems. Traditionally, model reduction techniques solve this problem by deriving low-dimensional models on the basis of the given high-dimensional model through suitable operations such as projection. In the field of systems and control, roughly two classes of model reduction techniques can be distinguished for linear systems: methods based on energy functions such as balanced truncation \cite{moore1981principal,1102945,4047847,benner2011lyapunov} and optimal Hankel norm approximation  \cite{glover1984all}, and methods based on interpolation and/or moment matching  \cite{feldmann1995efficient,grimme1997krylov,gallivan2004model,astolfi2010model}, sometimes also referred to as Krylov methods. Extensions to nonlinear systems have emerged in the form of nonlinear balancing methods \cite{scherpen1993balancing,besselink2014Incremental} and nonlinear moment matching techniques \cite{astolfi2010model,ionescu2015nonlinear}.
We refer the reader to \cite{antoulas2005approximation,antoulasInterpolatorybook2020,benner2015survey} for details on a variety of existing model reduction  methods.

Since recently, the problem of \emph{data-driven} model reduction is attracting increasing attention, partly motivated by the widespread availability of measurement data. Here, low-order models are constructed directly on the basis of measurement data, thus not requiring the availability of a high-order model. We emphasize that these data-driven model reduction approaches differ from traditional approaches in which, first, a high-order model is derived using system identification techniques and, second, existing model-based techniques for model reduction are used. Nevertheless, standard (model-based) model reduction techniques have inspired various data-driven techniques.  

First, in the class of energy-based methods for linear systems, to which this paper belongs, 
 \cite{rapisarda2011identification,markovsky2005algorithms} propose a data-driven balanced truncation method from \emph{persistently exciting} data and 
 \cite{gosea2021data} estimates Gramians from frequency and time-domain data based on their quadrature form. For nonlinear systems, empirical balanced truncation is presented in \cite{lall2002subspace,kawano2021empirical}, whereas data-driven reduction for monotone nonlinear systems is considered in \cite{Kawano2019monotone}. Second, in the class of interpolatory methods, we begin by mentioning contributions to data-driven  reduction methods on the basis of frequency-domain data by the Loewner framework \cite{mayo2007framework}. In this method, noise-free frequency-domain data are formulated to construct Loewner matrix pencils to enable the construction of state-space models. 
 Extensions of this approach aim at constructing 
 reduced-order models preserving stability \cite{gosea2016stability} and achieving $\HHH_{2}$-optimality \cite{beattie2012realization}.   In addition, the use of Loewner methods based on time-domain data and noisy frequency-domain data is pursued in \cite{peherstorfer2017data} and \cite{drmavc2019learning}, respectively. Besides the Loewner framework, data-driven moment matching techniques have been presented in \cite{astolfi2010model} and \cite{burohman2020data}, where the latter exploits the so-called data informativity framework.
  
 Despite these developments, existing methods for data-driven model reduction do often not allow for guaranteeing system properties such as asymptotic stability and do not provide an error bound, especially when the available data is subject to noise. In this paper, we develop a data-driven reduction technique that provides such guarantees on the low-order model, even for noisy data. Specifically, this paper has the following contributions. 
 
 First, we introduce the concept of  \emph{data reduction} via a Petrov-Galerkin projection. Following the data informativity framework of \cite{van2019data}, we characterize the class of systems that are consistent with the measurement data for a given noise model in terms of a quadratic matrix inequality (QMI). Then, we define the class of reduced-order systems as the set of systems obtained by applying the Petrov-Galerkin projection to all systems explaining the noisy data. Importantly, we show that this class of reduced-order systems can again be characterized in terms of a quadratic matrix inequality, but one of \emph{lower dimension}. As the relevant matrix variables in this QMI depend only on the measurement data, noise model and projection matrices, this can be regarded as \emph{data reduction}.
 
 The second contribution of this paper is the development of a data-driven generalized balanced truncation method. Here, we characterize the set of all reduced-order models obtained by applying generalized balanced truncation to the class of systems explaining the noisy data. This relies on the following two steps. 
 As the first step, we give necessary and sufficient conditions for all systems explaining the data to have a \emph{common} generalized controllability and \emph{common} generalized observability Gramian. In this case, we say that the data are \emph{informative} for generalized Lyapunov balancing. These conditions heavily rely on the so-called matrix S-lemma from \cite{van2020noisy} and again build on the data informativity framework of \cite{van2019data}. We note that this framework has also been successfully applied in solving various control problems, e.g., data-driven $\HHH_{2}$ and $\HHH_\infty$ control \cite{van2020noisy}. The second step comprises the use of the common generalized Gramians to obtain the Petrov-Galerkin projection that achieves (generalized) balanced truncation (see \cite{dullerud2013course} for details on generalized balanced truncation for a single given system). This allows for the application of the data reduction concept and yields the desired class of reduced-order models in terms of a low-dimensional quadratic matrix inequality.

  This data-driven model reduction procedure has various desirable properties by virtue of the inherent advantages of using a balancing-type reduction method. Namely, all reduced-order models are guaranteed to be asymptotically stable and satisfy an \emph{a priori} error bound.
  However, the ordinary a priori upper bound from model-based reduction methods, e.g., \cite{antoulas2005approximation,dullerud2013course}, does not determine the error between a selected reduced-order model (from the class of reduced-order models) to the true system generating the data because the true system is unknown. Therefore, as the final contribution of this paper, we provide two alternative error bounds. First, we compute a uniform a priori upper bound, i.e., an error bound that holds for any chosen high-order system explaining the data and any reduced-order model. The computation of this error bound again exploits the QMI characterization of the class of (reduced-order) systems, together with the bounded real lemma. Second, we also present an a posteriori error bound that is uniform over all systems explaining the data for a \emph{selected} reduced-order system.
 
 The remainder of this paper is organized as follows: In Section~\ref{sec:prel}, we provide appropriate background material on Petrov-Galerkin model reduction and generalized balanced truncation. Section~ \ref{Sec:DDpetrovgalerkin} deals with the data reduction problem through the Petrov-Galerkin projection. The problem formulation of informativity for  Lyapunov balancing is given in Section~\ref{sec:problemGLB}, followed by the main results containing necessary and sufficient conditions for generalized Lyapunov balancing, a characterization of the set of reduced-order models and error bounds in Sections~\ref{section:main}, \ref{section:set_red}, and \ref{section:errorbound}, respectively. In Section~\ref{sec:example}, an illustrative example is provided to show how the set of reduced-order models is extracted from data. Finally, the paper closes with some concluding remarks in Section~\ref{sec:conclusion}. For the sake of completeness, some important results and proofs are presented in the Appendix.
 	
\textit{Notation}. We denote $M>0$  $(M \geq 0)$ and $M<0$  $(M \leq 0)$ for positive and negative ({semi-}) definite symmetric matrices, respectively. We denote the number of negative, zero, and positive eigenvalues of  a symmetric matrix $M$ by $\upsilon_-(M)$, $\upsilon_0(M)$, and $\upsilon_+(M)$, respectively. The {\em inertia\/} of $M$ is denoted by $\In(M)=(\upsilon_-(M),\upsilon_0(M),\upsilon_+(M))$. 
For a symmetric matrix $M$ partitioned as
$$
M=\bbm M_{11} & M_{12} \\ M_{12}^\top & M_{22} \ebm,
$$
its Schur complement with respect to $M_{22}$ is denoted by $M | M_{22}$, i.e., $M | M_{22}:=M_{11}-M_{12}M_{22}^{-1}M_{12}^\top$.
For a square matrix $A$, its spectral radius is denoted by $\rho(A)$ and  the sum of its main diagonal elements is denoted by $\mathrm{trace}(A)$. We denote $\blkdiag(A_1,A_2, \ldots,A_n)$ for a block 
diagonal matrix of the form
$$
\bbm A_1 & 0 & \cdots  & 0 \\ 0 & A_2 & \cdots & 0 \\ \vdots & \vdots & \ddots & \vdots  \\ 0 & 0 &  \cdots & A_n \ebm.
$$
The matrix $I_j$ denotes the identity matrix of size $j$. 
For a discrete-time linear system $\Sigma$, its $\HHH_{\infty}$-norm is denoted by $\lVert \Sigma \rVert_{\HHH_{\infty}}$. For a system $\Sigma$ having realization $(A,B,C,D)$ and transfer function $G(z)=C(zI-A)^{-1}B+D$, the norm $\lVert \Sigma \rVert_{\HHH_{\infty}}$ is defined by $\lVert \Sigma \rVert_{\HHH_{\infty}}=\sup_{\omega \in \R} \lVert G(e^{i\omega}) \rVert$. 

\section{Preliminaries} \label{sec:prel}
\subsection{Model reduction via a Petrov-Galerkin projection} \label{sec:Pre sub:mor}
Consider the discrete-time input/state/output system
\begin{align}\label{e:sys4proj}
\Sigma: \quad \begin{split}
\bm x(k+1)&=A \bm x(k)+B \bm u(k),\\
\bm y(k)&=C \bm x(k)+D \bm u(k),
\end{split}
\end{align}
with input $\bm u \in \R^m$, state $\bm x \in \R^n$ and output $\bm y \in\R^{p}$. Let $\hatW,\hatV \in \R^{n\times r}$ be matrices such that $\hatW^\top \hatV =I$ and $r<n$. A reduced-order model obtained via a Petrov-Galerkin projection is given by
\begin{align}\label{e:redfromproj}
\hat{\Sigma}: \quad \begin{split}
\bm \hatx(k+1)&=\hatW^\top A \hatV \bm \hatx(k)+\hatW^\top B \bm u(k),\\
\bm y(k)&=C \hatV\bm \hatx(k)+D \bm u(k)
\end{split}
\end{align}
where $\bm \hatx \in \R^r$ denotes the state of the reduced-order model. The Petrov-Galerkin projection method provides a general framework for model order reduction. Namely, many reduction techniques, including Gramian- and Krylov-based methods can be regarded as Petrov-Galerkin projections with appropriate choice of $\hatV$ and $\hatW$, see, e.g., \cite{antoulas2005approximation}.

\subsection{Generalized Lyapunov balancing }\label{section:GLB}
Lyapunov balancing is a popular method for model reduction, which is slightly generalized in the concept of generalized Lyapunov balancing (GLB). We begin our discussion on GLB by introducing the notion of a generalized Gramian. Consider the discrete-time system \eqref{e:sys4proj}. Then, a matrix  $P=P^\top >0$ satisfying
\begin{equation}\nonumber
APA^\top -P+BB^\top  < 0
\end{equation}
and $Q=Q^\top >0$ satisfying 
\begin{equation}\nonumber
A^\top QA-Q+C^\top C< 0,
\end{equation}
are called the \emph{generalized controllability Gramian}  and \emph{generalized observability Gramian}, respectively.  This is a strict version of the definition of generalized Gramians given in \cite{dullerud2013course}. The generalized Gramians are lower bounded by the ordinary Gramians, i.e., $P> P_0$ and $Q> Q_0$, where $Q_0$ and $P_0$ are the solutions of the corresponding Lyapunov \emph{equations}. Note that asymptotic stability of $\Sigma$, i.e., $\rho(A) < 1$, where $\rho(A)$ represents the spectral radius of $A$, is a necessary condition for the generalized Gramians to exist.

These generalized Gramians can be used to obtain a so-called balanced realization, that is, a realization of \eqref{e:sys4proj} for which the corresponding generalized Gramians are equal and diagonal. Specifically, by \cite[Lemma~7.3]{antoulas2005approximation}, there exists a nonsingular matrix $T$ such that $T P T^\top = T^{-\top} Q T^{-1} = \Sigma_H$ where $\Sigma_H$ is a diagonal matrix of the \emph{generalized} Hankel singular values of  $\Sigma$ in \eqref{e:sys4proj}, i.e.,
\beq \label{e:Hankelsingval0}
\Sigma_H := \blkdiag ( \sigma_1I_{m_1},
\sigma_2 I_{m_2},
\ldots,
\sigma_\kappa I_{m_\kappa}
),
\eeq
with 
$
\sigma_1 > \sigma_2 > \cdots > \sigma_\kappa >0, 
$ and where $m_i$ denotes the multiplicity of $\sigma_i$ for $i=1,\dots,\kappa$ satisfying $n=\sum_{i=1}^{\kappa} m_i$.
In this case, we say that such a realization is \emph{balanced in the sense of generalized Lyapunov balancing}. 
In particular, the balanced realization is given by
\beq \label{e:ss_bal0}
	\! A_{\mathrm{bal}}\!:=T AT^{-1}\!,\ B_{\mathrm{bal}} \! :=T B,\
	C_{\mathrm{bal}}\! :=CT^{-1}, \ D_{\mathrm{bal}}\! :=\! D. \! \!
\eeq


Finally, the reduced-order model via balanced truncation is obtained by truncating the balanced system $(A_\mathrm{bal},B_\mathrm{bal},C_\mathrm{bal},D_\mathrm{bal})$. Indeed, balanced-truncation  model reduction is essentially a Petrov-Galerkin projection.
Namely, after introducing the matrix $\Pi \in \R^{n \times r}$ as 
\beq \nonumber 
\Pi := \bbm I_r \\ 0 \ebm,
\eeq 
the projection matrices $\hatV=T^{-1}\Pi$ and $\hatW=T^\top\Pi$ satisfy $\hatW^\top \hatV=I$ and the reduced-order model \eqref{e:redfromproj} is equal to the model obtained by balanced truncation (to order $r$).


Generalized balanced truncation guarantees the preservation of some relevant system properties, similarly as in ordinary balanced truncation, see \cite[Prop. 4.19]{dullerud2013course} for continuous-time systems. The discrete-time version is presented without proof below.
\bprop \label{prop:gen_bal}
Consider  the system $\Sigma$ given in \eqref{e:sys4proj}.  Let $\hat{\Sigma}$ of the form \eqref{e:redfromproj} be a reduced-order system of $\Sigma$ via generalized balanced truncation. Suppose that $\hat{\Sigma}$ is of order $r<n$ where $r = \sum_{i=1}^{\ell}m_i$ with $\ell < \kappa$. Then it is balanced in the  sense of GLB with $\rho(\hatW^\top A \hatV)  <1$ and 
\beq \nonumber 
\lVert \Sigma -\hat{\Sigma} \rVert_{\HHH_{\infty}} \leq 2 \sum_{i=\ell +1}^{\kappa}\sigma_i,
\eeq
where the $\sigma_i$'s are the neglected generalized Hankel singular values given in \eqref{e:Hankelsingval0}.
\eprop
\begin{figure*}[b]
	\par\noindent\rule{\textwidth}{0.5pt}
	\centering
	\newcounter{mytempeqncntR}
	\normalsize
	\setcounter{mytempeqncntR}{\value{equation}}
	\begin{equation}\label{e:NVW} 
	\tag{13} 
	N_{V,W} :=  \begin{pmat}[{|}]
	W^\top (N | N_{22} +N_{12}N_{22}^{-1} V (V^\top N_{22}^{-1}V)^{-1}V^\top N_{22}^{-1}N_{12}^\top )W &  W^\top N_{12} N_{22}^{-1} V (V^\top N_{22}^{-1}V)^{-1} \cr\-
	(V^\top N_{22}^{-1}V)^{-1}V^\top N_{22}^{-1}N_{12}^\top W  & (V^\top N_{22}^{-1} V)^{-1} \cr
	\end{pmat}
	\end{equation}
\end{figure*} 

\section{Data-driven Petrov-Galerkin projection}\label{Sec:DDpetrovgalerkin}
Consider the linear discrete-time input/state/output system
\begin{align}\label{e:tru-sys}
\Sigma_{\mathrm{true}}: \quad \begin{split}
\bm x(k+1)&=A_{\mathrm{true}} \bm x(k)+B_{\mathrm{true}} \bm u(k)+\bm w(k),\\
\bm y(k)&=C_{\mathrm{true}} \bm x(k)+D_{\mathrm{true}} \bm u(k)+\bm z(k),
\end{split}
\end{align}
where $(\bm u,\bm x,\bm y)\in\R^{m+n+p}$ are the input/state/output and $(\bm w,\bm z)\in\R^{n+p}$ are noise terms. Throughout the paper, we assume that the system matrices $(A_{\mathrm{true}}, B_{\mathrm{true}},C_{\mathrm{true}},D_{\mathrm{true}})$ and the noise $(\bm w,\bm z)$ are {\em unknown\/}. What is known instead are a finite number of input/state/output measurements harvested from the true system \eqref{e:tru-sys}:
\begin{align} 
u(0)&,u(1), \ldots, u(L-1), \nonumber \\
x(0)&,x(1), \ldots, x(L),\nonumber \\
y(0)&,y(1), \ldots, y(L-1). \nonumber
\end{align}
We collect these data in the matrices
\begin{align}
X&:=\begin{bmatrix}
x(0)&x(1)&\cdots& x(L) 
\end{bmatrix}, \nonumber\\
X_-& :=\begin{bmatrix}
x(0)&x(1)&\cdots& x(L-1)
\end{bmatrix},\nonumber\\
X_+& :=\begin{bmatrix}
x(1)&x(2)&\cdots& x(L)
\end{bmatrix}, \nonumber\\
U_-& :=\begin{bmatrix}
u(0)&u(1)&\cdots& u(L-1)
\end{bmatrix}, \nonumber\\
Y_-& :=\begin{bmatrix}
y(0)&y(1)&\cdots& y(L-1) 
\end{bmatrix}. \nonumber
\end{align}
Now, we can define the set of all systems that {\em explain\/} the data as 
\beq \nonumber 
\sig :=\set{\left(A,B,C,D\right)}{\begin{bmatrix} X_+\\Y_- \end{bmatrix}- \begin{bmatrix} A&B\\
C&D\end{bmatrix}\begin{bmatrix}X_-\\ U_- \end{bmatrix}\in\calN},
\eeq 
where $\calN\subseteq\R^{(n+p)\times L}$ captures a {\em noise model\/} such that 
\beq  \label{e:true in exp}
(A_{\mathrm{true}},B_{\mathrm{true}},C_{\mathrm{true}},D_{\mathrm{true}})\in\sig.
\eeq
In this paper, we work with a noise model that is described by a quadratic matrix inequality as 
\begin{equation}\label{eq:N1}
\calN\!:=\!\set{\!Z\in\R^{(n+p)\times L\!}}{
{\bbm
I \\ Z^\top 
\ebm\!}
^\top \!\!
\bbm
\Phi_{11} & \Phi_{12}\\
\Phi_{12}^\top  & \Phi_{22}
\ebm
\!\!
\bbm
I \\ Z^\top  
\ebm
\!\geq 0
},
\end{equation}
where $\Phi_{11}=\Phi_{11}^\top \in\R^{(n+p)\times(n+p)}$, $\Phi_{12}\in\R^{(n+p)\times L}$, and $\Phi_{22}=\Phi_{22}^\top \in\R^{L\times L}$. 

Throughout the paper, we make the following blanket assumption on the set $\calN$.
\bas \label{as:boundednonempty}
The set $\calN$ is bounded and has nonempty interior.
\eas
As shown in \cite{vanDissipativity2022mag}, one can verify this assumption by using the following lemma.
\setlength\arraycolsep{2pt}
\blem \label{l:bounded_nonempty}
The set $\calN$ given by \eqref{eq:N1}
is bounded and has nonempty interior if and only if $\Phi_{22}<0$ and ${\Phi_{11}-\Phi_{12}\Phi_{22}\inv\Phi_{12}^\top >0}$.
\elem

 It is clear from the definition of $\sig$ and \eqref{eq:N1} that $(A,B,C,D) \in \sig$ if and only if the following quadratic matrix inequality (QMI) is satisfied
 \begin{equation}\label{e:Ineq_noise}
 \setlength\arraycolsep{2pt}
 \systr^\top N\systr\geq 0,
 \end{equation}
 where
 \beq \label{e: N}
 N :=\begin{bmatrix}
 	I & 0  & X_{+} \\0 & I & Y _{-}\\ 0 &0 & -X_{-} \\ 0 & 0 & -U_{-}
 \end{bmatrix}\begin{bmatrix}
 	\Phi_{11}  & \Phi_{12} \\ \Phi_{12}^\top  & \Phi_{22}
 \end{bmatrix} \begin{bmatrix}
 	I & 0  & X_{+} \\0 & I & Y _{-}\\ 0 &0 & -X_{-} \\ 0 & 0 & -U_{-}
 \end{bmatrix}^\top.
 \eeq

In characterizing the set of systems that explain the data, one may wonder whether the set $\sig$ is bounded and has nonempty interior.
The following proposition provides the required condition, which solely relies on the data.
\bprop \label{l:N nonsingular}  
The set  $\sig$ is bounded and has nonempty interior if and only if  
there exists $\barS \in \R^{(n+m)\times(n+p)}$ such that
\beq \label{e: slater on N}
\bbma I \\ \barS \ebma^{\top} N \bbma I \\ \barS \ebma > 0.
\eeq 
and $\bbma X_{-} \\ U_{-} \ebma$ has full row rank.
\eprop
\begin{proof}
	The proof is presented in Appendix~\ref{app:N_nonsingular}.
\end{proof}
The condition \eqref{e: slater on N} is referred as the \emph{generalized Slater condition}.

In this paper, we are not interested in the true system \eqref{e:tru-sys} per se. Instead, we would like to find a \emph{reduced-order approximation} of \eqref{e:tru-sys} directly on the basis of  the available data.

As a first step, we consider a Petrov-Galerkin projection as in Section~\ref{sec:prel} and assume that the projection matrices $\hatW$ and $\hatV$ satisfying $\hatW^\top \hatV = I$ are given. Then, the set of reduced-order models of all systems explaining the data is defined as
\beq \label{e:set_sigred} \nonumber 
\!\!\sigred \!:=\!\set{\!(\hatW^\top A \hatV , \hatW^\top B, C\hatV,D)}{(A,B,C,D) \! \in \!\sig }\!.
\eeq 

The first main result of this paper is that the set $\sigred$ can itself be represented as a QMI of a similar form as \eqref{e:Ineq_noise}. This is formalized in the following theorem, whose proof can be found in Appendix~\ref{app_proj}.
\bthe \label{thm:ROM in QMI}
Consider the set $\sig$ of systems explaining the data. Suppose that there exists $\barS$ such that \eqref{e: slater on N} holds and the matrix $\bbma X_{-} \\ U_{-} \ebma$ has full row rank. Let $\hatW,\hatV \in \R^{n\times r}$ be such that $\hatW^\top \hatV =I$. Then,
the set $\sigred$ of reduced-order models of $\sig$ using  projection matrices $\hatW, \hatV$ satisfies
\beq \nonumber 
\!\sigred \!=\! \set{\!(\hatA,\hatB,\hatC,\hatD)\!}{\!\bbm I & 0 \\ 0 & I \\ \hatA^\top & \hatC^\top \\ \hatB^\top & \hatD\!^\top \ebm^\top \! \! N_{V,W} \bbm I & 0 \\ 0 & I \\ \hatA^\top & \hatC^\top \\ \hatB^\top & \hatD^\top \ebm\!\geq 0\!}\!,\!
\eeq 
 where $N_{V,W}$ is given by \eqref{e:NVW} with
 $$
 W:= \bbm \hatW & \\ & I_p \ebm \text{ and } V:= \bbm \hatV & \\ & I_m \ebm.
 $$
\ethe 

 Theorem~\ref{thm:ROM in QMI} has a nice interpretation in terms of \emph{data reduction}. Namely, the matrix $N_{V,W}$ characterizing all reduced-order models depends only on the projection matrices $\hatV, \hatW$ and the original data matrix $N$. As such, $N_{V,W}$ is constructed from the data and noise model only. Importantly, $N_{V,W}$ has a lower dimension than $N$ and can thus be regarded as a reduced data matrix. Hence, we can characterize all reduced-order models by directly reducing the data matrix $N$ rather than reducing individual systems $(A,B,C,D) \in \sig$. It is also worth mentioning that 
 \beq \nonumber
 (\hatW^\top A_\mathrm{true} \hatV,\hatW^\top B_\mathrm{true},  C_\mathrm{true} \hatV, D_\mathrm{true}) 
 \eeq
 i.e., the reduced-order model of the true system, is in $\sigred$.


 In this section, we have characterized reduced-order approximations of systems explaining the collected data for \emph{given} projection matrices $\hatW$ and $\hatV$. In the next section, we will choose the projection matrices on the basis of the available data by following a generalized balancing framework.
%

\section{Data-driven generalized balanced truncation} \label{sec:problemGLB}
In this section, we will introduce the notion of informativity for generalized Lyapunov balancing (GLB). Moreover, we give necessary and sufficient  conditions for informativity for  GLB, followed by the set of reduced-order models obtained from data-driven GLB and their error-bounds.
\begin{figure*}[b]
	\par\noindent\rule{\textwidth}{0.5pt}
	\centering
	\newcounter{mytempeqncnt}
	\normalsize
	\setcounter{mytempeqncnt}{\value{equation}}
	\setcounter{equation}{0}
	\begin{equation}\label{e:LMIs} \tag{16}
	\begin{pmat}[{|}] 
	\blkdiag(K_{11},(\frac{1}{2} - \mu) I_p , -K_{11}, -\gamma^{-2} I_m) & \blkdiag(K_{12}, - \mu I_p , -K_{12}, -\gamma^{-2} I_m) \cr\-
	\blkdiag(K_{12}^\top, -\mu I_p , -K_{12}^\top, -\gamma^{-2} I_m) & \blkdiag(K_{22}, (\frac{1}{2} - \mu) I_p , -K_{22}, -\gamma^{-2} I_m) \cr
	\end{pmat}
	- \blkdiag(\delta {N},\eta  {N}_{V,W} ) >0
	\end{equation}
	\setcounter{equation}{\value{mytempeqncnt}}
\end{figure*}

\subsection{Data informativity for generalized Lyapunov balancing}\label{section:main}
Based on Section \ref{section:GLB}, one can introduce the notion of informativity for GLB as follows.
\begin{definition}\label{def:Lyap_balancing}
We say that	the data $(U_{-},X,Y_-)$ are {\em informative for generalized Lyapunov balancing (GLB)\/} if there exist $P=P^\top >0$ and $Q=Q^\top >0$ such that
	\begin{equation}\label{e:Ineq_ctrbGramian}
	APA^\top -P+ BB^\top  < 0
	\end{equation}
	and
	\begin{equation}\label{e:Ineq_obsvGramian}
    A^\top QA-Q+ C^\top C < 0
	\end{equation}
	for all $\left(A,B,C,D\right)\in\sig$.
\end{definition}

\bre
It is known, see, e.g., \cite[Chapter~7]{antoulas2005approximation}, that the satisfaction of \eqref{e:Ineq_ctrbGramian} or \eqref{e:Ineq_obsvGramian} implies that all systems in $\sig$ are asymptotically stable.
\ere

From Definition~\ref{def:Lyap_balancing}, $P$ and $Q$ can be regarded as \emph{common} generalized controllability and observability Gramian, respectively, for all systems explaining the data. We thus formalize the following informativity and model reduction problems.
\begin{problem}
	 Find necessary and sufficient conditions under which the data $(U_{-},X,Y_{-})$ are informative for generalized Lyapunov balancing (GLB). If the data are informative for GLB, then characterize the reduced-order models via data-driven balanced truncation and provide error bounds with respect to the true system.
\end{problem}

Observe that the data are informative for GLB if and only if  QMI \eqref{e:Ineq_noise} implies the existence of positive definite matrices $P$ and $Q$ such that \eqref{e:Ineq_ctrbGramian} and \eqref{e:Ineq_obsvGramian} hold.  Such QMI implications can be viewed as a generalization of the classical S-lemma \cite{yakubovich1977s} and have been investigated in \cite{van2020noisy}. Based on the results of \cite{van2020noisy} and \cite{vanDissipativity2022mag}  (see Appendix~\ref{app:Slemma}), data informativity for GLB can be fully characterized in terms of feasibility of certain LMIs as stated next. 
\bthe \label{thm:inf_bal}
Suppose that  there exists $\barS$ such that \eqref{e: slater on N} holds. Define
\beq \nonumber
N_{\calC}:= \bbm I_n & 0 \\ 0 & 0 \\ 0 & I_{n+m} \ebm^\top N \bbm I_n & 0 \\ 0 & 0 \\ 0 & I_{n+m} \ebm
\eeq 
and
\beq \nonumber 
{N}_{\calO} :=\bbm I_n & 0 \\ 0 & 0 \\ 0 & I_{n+p} \ebm^\top N^{\sharp} \bbm I_n & 0 \\ 0 & 0 \\ 0 & I_{n+p} \ebm,
\eeq
where 
\beq \nonumber 
N^{\sharp} :=\begin{bmatrix} 0 & -I_{n+m} \\ I_{n+p} & 0 \end{bmatrix} N^{-1} \begin{bmatrix} 0 & -I_{n+p} \\ I_{n+m} & 0 \end{bmatrix}.
\eeq 
Then, the data $(U_-,X,Y_-)$ are informative for generalized Lyapunov balancing if and only if
\newcounter{counterafterNVW}
\normalsize
\setcounter{counterafterNVW}{\value{equation}}
	\setcounter{equation}{13}
\begin{enumerate}
	\renewcommand{\theenumi}{(\roman{enumi})}
	\renewcommand{\labelenumi}{\theenumi}
	\item \label{mainthm1} $\bbma X_- \\ U_- \ebma$ has full row rank,
	\item \label{mainthm2} there exists ${P}={P}^\top >0$ and a scalar $\alpha > 0$ such that
 	\begin{equation}\label{e:Slemma_ctrb}
 	\begin{pmat}[{..}]
 	{P} & 0 & 0 \cr
 	0 &  -P & 0 \cr
 	0 &  0 & -I_m  \cr
 	\end{pmat}-\alpha  {N}_\calC > 0,
 	\end{equation}
 	\item \label{mainthm3}there exists ${Q}={Q}^\top >0$ and a scalar $\beta > 0$ such that
 		\begin{equation}\label{e:Slemma_obsv}
 	\begin{pmat}[{..}]
 	{Q} & 0 & 0 \cr
 	0 & -{Q} & 0 \cr
 	0  & 0 & -I_p  \cr
 	\end{pmat}-\beta  N_{\calO} > 0.
 	\end{equation}
\end{enumerate}
\ethe
\begin{proof}
The proof can be found in Appendix~\ref{app:thm_inf}.
\end{proof}

A direct consequence of data informativity for GLB is that all systems explaining the data have \emph{common} generalized Gramians $P$ and $Q$. As a result, all systems in $\sig$ are balanced by a common balancing transformation matrix $T$ satisfying $T P T^\top = T^{-\top} Q T^{-1} = \Sigma_H$ where  $\Sigma_H$ is a matrix of the form \eqref{e:Hankelsingval0}, containing the common generalized Hankel singular values. Next, we note that the balanced realizations of all systems explaining the data can be constructed from \eqref{e:ss_bal0} where $(A,B,C,D) \in \sig$. 

\bre
Since $P$ and $Q$ satisfying \eqref{e:Slemma_ctrb} and \eqref{e:Slemma_obsv}, respectively, are lower bounded by the ordinary Gramians of 
the true system, (see the discussion in Section~\ref{section:GLB}), smaller $P$ and $Q$ are expected to yield a balancing that is ``closer'' to the ordinary balancing of the true system. Therefore, one may solve LMIs \eqref{e:Slemma_ctrb} and \eqref{e:Slemma_obsv} by minimizing $\mathrm{trace}(P)$ and $\mathrm{trace}(Q)$ to expect a better reduced-order approximation.
\ere

In the next section, we will use these balanced realizations to obtain  reduced-order models directly from data. 

\subsection{Reduced-order models} \label{section:set_red}
By applying the Petrov-Galerkin projection, the reduced-order models of all  systems in $\sig$ via generalized balanced truncation
are contained in the set
\beq \label{e:set_sigredBT} \nonumber 
\sigredbt :=\set{\!(\hatW^\top A \hatV , \hatW^\top B, C\hatV,D)}{(A,B,C,D) \! \in \!\sig }\!
\eeq 
where $\hatV=T^{-1}\Pi$ and $\hatW=T^\top\Pi$  with $\Pi$ given by 
\beq \nonumber
\Pi := \bbm I_r \\ 0 \ebm
\eeq 
and $T$ is obtained from the common generalized Gramians $P$ and $Q$ for all systems in $\sig$.
Based on Theorem~\ref{thm:ROM in QMI}, we can characterize the set  $\sigredbt$ in terms of a quadratic matrix inequality. We formalize this fact in the following corollary of Theorem~\ref{thm:ROM in QMI}. 
\begin{corollary}\label{cor:red_in_QMI}
Suppose that there exists $\barS$ such that \eqref{e: slater on N} holds and the data $(U_-,X,Y_-)$ are informative for generalized Lyapunov balancing with $T$ the corresponding balancing transformation. Then, 
\beq \nonumber 
\!\sigredbt = \set{\!(\hatA,\hatB,\hatC,\hatD)\!}{\! \systrred^\top \! \! \! \! N_{V,W} \systrred\!\geq 0}\!,
\eeq 
where $N_{V,W}$ is given by \eqref{e:NVW} with 
\beq \nonumber 
W=\bbm T^\top \Pi  & \\ & I_p \ebm, \  V=\bbm T^{-1}\Pi  & \\ & I _m\ebm  
\text{ and } \Pi := \bbm I_r \\ 0 \ebm .
\eeq 
\end{corollary}

\begin{figure*}[b]
	\par\noindent\rule{\textwidth}{0.5pt}
	\centering
	\newcounter{mytempeqncnt2}
	\normalsize
	\setcounter{mytempeqncnt2}{\value{equation}}
	\setcounter{equation}{0}
	\begin{equation}\label{e:LMI_0} \tag{21}
	\begin{pmat}[{...|.}]
	K_{11} & 0  & 0  &  0 & K_{12}   \cr
	0& I_p-\hatC_0 K_{22} \hatC_0^\top - \gamma_0^{-2}\hatD_0\hatD_0^\top  & \hatC_0 K_{12}^\top & \gamma_0^{-2} \hatD_0 &  \hatC_0 K_{22}\hatA_0^\top + \gamma_0^{-2} \hatD_0\hatB_0^\top  \cr
	0 & K_{12}\hatC_0^\top & -K_{11} & 0 &  -K_{12} \hatA_0^\top \cr
	0 & \gamma_0^{-2}\hatD_0^\top  & 0  & -\gamma_0^{-2} I_m &  -\gamma_0^{-2}\hatB_0^\top \cr\-
	K_{12}^\top & \hatA_0 K_{22}\hatC_0^\top + \gamma_0^{-2}\hatB_0 \hatD_0^\top  & -\hatA_0 K_{12}^\top & -\gamma_0^{-2} \hatB_0  &   K_{22} - \hatA_0 K_{22} \hatA_0^\top - \gamma_0^{-2} \hatB_0 \hatB_0^\top \cr
	\end{pmat} - 
	\begin{pmat}[{|}]
	\delta {N} &  0\cr\-
	0 & 0  \cr
	\end{pmat}  >0
	\end{equation}
	\setcounter{equation}{\value{mytempeqncnt2}}
\end{figure*}

Now, we can see that the set $\sigredbt$ characterizes a \emph{data reduction} by generalized balanced truncation. Namely, the reduced matrix $N_{V,W}$ depends only on the data matrix $N$ and projection matrices $\hatV=T^{-1}\Pi$ and $\hatW=T^\top\Pi$, where now these projection matrices are also derived from the data only via Theorem~\ref{thm:inf_bal}.


From the definition of $\sigredbt$ above, suppose that $(\hatA,\hatB,\hatC,\hatD)  \in \sigredbt$, then it is always a truncation of a model in $\sig$ by generalized balanced truncation.  Therefore, any $(\hatA,\hatB,\hatC,\hatD)  \in \sigredbt$  satisfies the guaranteed properties in Proposition~\ref{prop:gen_bal}. Namely,  $ \rho(\hatA)  <1$ and the $\HHH_{\infty}$-norm  error between a system in $\sigredbt$ and its \emph{corresponding} system in $\sig$ is upper bounded by the neglected common generalized Hankel singular values, i.e.,
\newcounter{counterafterLMI1}
\normalsize
\setcounter{counterafterLMI1}{\value{equation}}
\setcounter{equation}{16}
\beq \label{e:upper_bound_bal_data}
\lVert \Sigma -\hat{\Sigma} \rVert_{\HHH_{\infty}} \leq 2 \sum_{i=\ell +1}^{\kappa}\sigma_i
\eeq
where $\Sigma$ and $\hat{\Sigma}$ are systems whose realizations are in $\sig$ and $\sigredbt$, respectively, and $\sigma_i$ is defined similarly as in Proposition~\ref{prop:gen_bal}.

However, the bound \eqref{e:upper_bound_bal_data} has little practical relevance as it characterizes the error between \emph{one} reduced-order system in $\sigredbt$ and its \emph{corresponding} high-order system in $\sig$.  Instead, recall that we are interested in a reduced-order approximation of the true system \eqref{e:tru-sys}. This system is unknown, but is guaranteed to satisfy \eqref{e:true in exp}. As also the corresponding reduced-order system is unknown (but in $\sigredbt$), a practical relevant error bound should hold for \emph{any} selection of a high-order system (from $\sig$) and \emph{any} reduced-order system (from $\sigredbt$). The following section provides such bounds.
%

\subsection{Distance to the true system} \label{section:errorbound}
Suppose that we consider a reduced-order system of  order $r<n$ given by  $\hat{\Sigma}_0$ with realization $(\hatA_0,\hatB_0,\hatC_0,\hatD_0) \in \sigredbt$.
We will use $\hat{\Sigma}_0$ as an approximation of  the unknown true system $\Sigma_{\mathrm{true}}$. To evaluate the quality of this approximation, note that
\beq \label{e:sup} \lVert \hat{\Sigma}_0- {\Sigma}_{\mathrm{true}} \rVert_{\HHH_\infty}\leq   \sup \big\{ \| \hat{\Sigma}- \Sigma \|_{\HHH_\infty} \! :\! {\Sigma} \in \sig, \hat{\Sigma} \in \sigredbt \big \}.\eeq 
Here, we have used the small abuse of notation ${\Sigma} \in \sig$ to mean $(A,B,C,D) \in \sig$, where $(A,B,C,D)$ is a realization of $\Sigma$.

In this section, we aim to compute a bound on the right-hand side of  \eqref{e:sup} on the basis of the available data only. The computation of this is stated in the following result.
\bthe \label{thm:maxG1G2}
The bound  
\beq \nonumber 
\| \hat{\Sigma}- \Sigma \|_{\HHH_\infty}< \gamma
\eeq 
holds  for any $\Sigma \in \sig$ and any $\hat{\Sigma} \in \sigredbt$  
  if and only if there exist a  matrix $K=K^\top >0$ in $\R^{(n+r) \times (n+r)}$ partitioned as 
\beq \nonumber
K=\bbm K_{11} & K_{12} \\ K_{12}^\top & K_{22} \ebm, \text{ with } K_{11} \in \R^{n \times n},
\eeq
and scalars  $\delta>0$, $\eta>0$ and $\mu$ such that  \eqref{e:LMIs} holds, where $N$ and $N_{V,W}$ are given by \eqref{e: N} and \eqref{e:NVW}, respectively.
\ethe
\begin{proof}
	The proof can be found in Appendix~\ref{app:proofbound}.
\end{proof}


In order to obtain the smallest upper bound, i.e., the smallest $\gamma$ such that the conditions in Theorem~\ref{thm:maxG1G2} hold, 
one may solve the following  semidefinite program \cite[Sect. 6.4]{semidefiniteprogYuriiArkadii}:
\begin{subequations}\label{prob:opt}
	\begin{alignat}{2}
	&\!\min_{K,\delta,\eta,\mu} \quad \gamma  \\
	& \text{subject to} & \quad K >0, \delta > 0, \eta > 0 \text{ and }\eqref{e:LMIs}.
	\end{alignat}
\end{subequations}
In conclusion, the solution of  \eqref{prob:opt} gives $\lVert \hat{\Sigma}_0 - \Sigma_{\mathrm{true}} \rVert < \gamma$. Note that this upper bound is uniform for any $\hat{\Sigma}_0$  picked from $\sigredbt$. Therefore, it can be regarded as an \emph{a priori} error bound.

\subsection*{A posteriori error bound}
Now suppose we pick a \emph{known} system $\hat{\Sigma}_0$ from  $\sigredbt$. Let its realization be $(\hatA_0,\hatB_0,\hatC_0,\hatD_0)$. 
An a posteriori error bound is computed to measure the error between this known system $\hat{\Sigma}_0$ and the true system $\Sigma_{\mathrm{true}}$. However, since the true system is unknown, this error cannot be directly computed. 

Fortunately, we know that 
$$\lVert \hat{\Sigma}_0 - \Sigma_\mathrm{true} \rVert_{\HHH_{\infty}} \leq \sup \set{\lVert \hat{\Sigma}_0 - \Sigma \rVert_{\HHH_{\infty}} }{\Sigma \in \sig}.$$
The following proposition gives the computation of an upper bound of $\lVert \hat{\Sigma}_0 - \Sigma \rVert_{\HHH_{\infty}}$ for any $\Sigma \in \sig$.

\bprop \label{prop:upper_boundG0}
Let $\hat{\Sigma}_0$ be a given reduced-order model with realization $(\hatA_0,\hatB_0,\hatC_0,\hatD_0) \in \sigredbt$. Then,  the bound 
\beq \nonumber 
\lVert \hat{\Sigma}_0 -\Sigma \rVert_{\HHH_{\infty}} < \gamma_0
\eeq 
 holds for any $\Sigma \in \sig$  if and only if there exist ${K}={K}^\top >0$ in $\R^{(n+r) \times (n+r)}$ partitioned as 
\beq \nonumber
K=\bbm K_{11} & K_{12} \\ K_{12}^\top & K_{22} \ebm, \text{ with } K_{11} \in \R^{n \times n},
\eeq
 and a scalar ${\delta} > 0$ such that \eqref{e:LMI_0} holds.
\eprop
\begin{proof}
	The proof is presented in Appendix~\ref{app:boundsys0}.
\end{proof}

Similar to before, one may obtain the smallest upper bound by solving the following semidefinite program:
\begin{subequations}\label{prob:opt0}
	\begin{alignat}{2}
	&\!\min_{{K},{\delta}}  \ \gamma_0  \\
	& \text{subject to} & \quad {K} >0, {\delta} > 0, \text{ and }\eqref{e:LMI_0}.
	\end{alignat}
\end{subequations}
Finally, we have $\lVert \hat{\Sigma}_0 - \Sigma_{\mathrm{true}} \rVert_{\HHH_{\infty}} < \gamma_0$. Note that this a posteriori error bound  holds for a specific $\hat{\Sigma}_0$. 
It follows readily that $\gamma_0 \leq \gamma$ where $\gamma$ and $\gamma_0$ are the solutions of \eqref{prob:opt} and \eqref{prob:opt0}, respectively, since the bound $\gamma_0$ holds for a specific $\hat{\Sigma}_0 \in \sigredbt$ while the bound $\gamma$ holds for any $\hat{\Sigma}$ taken from $\sigredbt$.
%
%
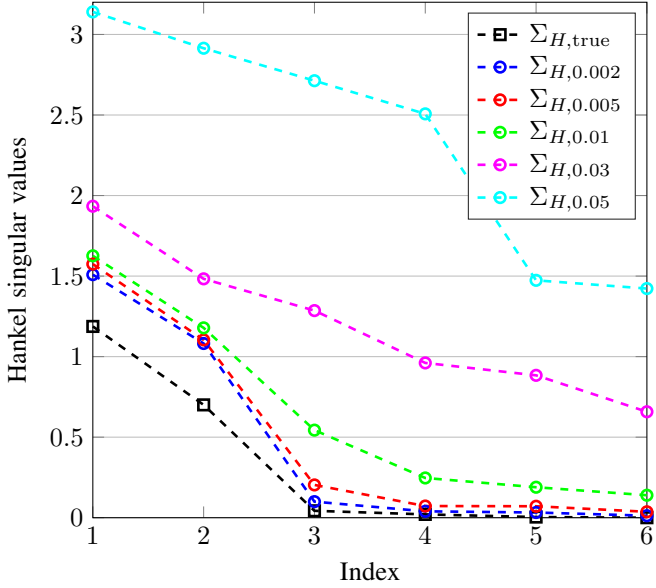
\begin{figure}[t]
%
%
\definecolor{mycolor1}{rgb}{1.00000,0.00000,1.00000}%
\definecolor{mycolor2}{rgb}{0.00000,1.00000,1.00000}%
\begin{tikzpicture}

\begin{axis}[%
width=2.9in,
height=2.7in,
at={(0.011in,0.25in)},
scale only axis,
xmin=1,
xmax=6,
xlabel={Index},
ymin=0,
ymax=3.2,
ylabel={Hankel singular values},
ymajorgrids,
axis background/.style={fill=white},
legend style={legend cell align=left, align=left, draw=white!15!black}
]
\addplot [color=black, dashed, line width=1.0pt, mark=square, mark options={solid, black}]
  table[row sep=crcr]{%
1	1.18744248023713\\
2	0.700455898716033\\
3	0.0428440944587903\\
4	0.0205769916903347\\
5	0.00411116520165414\\
6	3.79308258700447e-05\\
};
\addlegendentry{$\Sigma_{H,\mathrm{true}}$}

\addplot [color=blue,dashed, line width=1.0pt, mark=o, mark options={solid, blue}]
  table[row sep=crcr]{%
1	1.50889372907176\\
2	1.08156735165088\\
3	0.0996667796145758\\
4	0.0396314056657524\\
5	0.0322409995875914\\
6	0.0125057295553518\\
};
\addlegendentry{$\Sigma_{H,0.002}$}

\addplot [color=red, dashed,line width=1.0pt, mark=o, mark options={solid, red}]
  table[row sep=crcr]{%
1	1.57447623296142\\
2	1.10378990781595\\
3	0.203482153544972\\
4	0.0724532507735195\\
5	0.0710935344720486\\
6	0.0367128831821971\\
};
\addlegendentry{$\Sigma_{H,0.005}$}

\addplot [color=green,dashed, line width=1.0pt, mark=o, mark options={solid, green}]
  table[row sep=crcr]{%
1	1.62593289540247\\
2	1.17853649486511\\
3	0.543629119713791\\
4	0.246724042519139\\
5	0.189178255205048\\
6	0.139480180303037\\
};
\addlegendentry{$\Sigma_{H,0.01}$}

\addplot [color=mycolor1, dashed,line width=1.0pt, mark=o, mark options={solid, mycolor1}]
  table[row sep=crcr]{%
1	1.93367123241255\\
2	1.48294820731704\\
3	1.28578267800964\\
4	0.960901987361825\\
5	0.883861596096025\\
6	0.657659914249099\\
};
\addlegendentry{$\Sigma_{H,0.03}$}

\addplot [color=mycolor2, dashed, line width=1.0pt, mark=o, mark options={solid, mycolor2}]
  table[row sep=crcr]{%
1	3.1400441396493\\
2	2.91446182562371\\
3	2.71293627481078\\
4	2.50707400945421\\
5	1.47366254168939\\
6	1.42257331532415\\
};
\addlegendentry{$\Sigma_{H,0.05}$}

\end{axis}

\begin{axis}[%
width=2.778in,
height=2.833in,
at={(0in,0in)},
scale only axis,
xmin=0,
xmax=1,
ymin=0,
ymax=1,
axis line style={draw=none},
ticks=none,
axis x line*=bottom,
axis y line*=left,
legend style={legend cell align=left, align=left, draw=white!15!black}
]
\end{axis}
\end{tikzpicture}%
	\caption{The Hankel singular values of the true system denoted by $\Sigma_{H,\mathrm{true}}$ and generalized Hankel singular values of all systems explaining the data for different noise levels denoted by $\Sigma_{H,\sigma}$.  The values $\Sigma_{H,\mathrm{true}}$ are computed via the ordinary balancing procedure, i.e., balancing by Lyapunov equations, and using the system matrices of $\Sigma_{\mathrm{true}}$. }
	\label{fig:Hankel}
\end{figure}
\section{Illustrative example}\label{sec:example}
Consider a continuous-time system of  a cart with double pendulum, see \cite{ionescu2015two} for details. After discretizing this model using the zero-order hold method with sampling time $0.5$ seconds, we obtain  a true discrete-time system of the form \eqref{e:tru-sys} where  $A_{\mathrm{true}}$, $B_{\mathrm{true}}$, $C_{\mathrm{true}}$, and $D_{\mathrm{true}}$ are given by
\beq \nonumber 
\bbm 
0.9299  &  0.4160 &   0.7447   & 0.2291  &  0.2452  &  0.0592 \\
-0.1869  &  0.7430   & 0.3318 &   0.7617 &   1.0859 &   0.3560 \\
0.0380 &   0.0477 &  -0.3644 &   0.0647  &  0.1370  &  0.0766 \\
0.0169  &  0.0549 &  -0.0972  & -0.3693  &  -0.8685  &  0.0484 \\
0.0250 &   0.0285  &  0.2741  &  0.1393 &  -0.0474   & 0.1615 \\
0.1108  &  0.1358  & -1.7370   & 0.1855 &  -1.8002  & -0.2311
\ebm,
\eeq 
\beq \nonumber \bbm 
0.0701 &
0.1869 &
-0.0380 &
-0.0169 &
-0.0250 &
-0.1108 
\ebm^\top, \eeq 
$\bbm 1 & 0 & 0 & 0 & 0 & 0 \ebm$, and $0$, respectively.
%
	\newcounter{counterafterLMI2}
\normalsize
\setcounter{counterafterLMI2}{\value{equation}}
\setcounter{equation}{21}

To illustrate the data-driven model reduction from noisy data, we collect input/state/output data of system \eqref{e:tru-sys} up to $L=200$ for  input signal
\beq \label{e:input_signal}
u(k)=2\sin(k)+\cos(0.5 k)
\eeq 
and a random initial condition $x(0)$ which follows a Gaussian distribution with zero mean and unit variance.
In addition, we take the noise $\bm w$ and $\bm z$ in \eqref{e:tru-sys} to be Gaussian with zero mean and variance $\sigma^2$. For a realization of this noise, we obtain the data matrices $U_-, X_-, X_+, Y_-$. 
In the remainder of this example, we assume that the noise samples satisfy noise model \eqref{eq:N1} with $\Phi_{11} = 1.35 \sigma^2 I$, $\Phi=0$ and $\Phi_{22} = -I$.
We simulate the noise  with different levels: $\sigma \in \{0.002,0.005,0.01,0.03,0.05\}$ and it has been checked that they satisfy the noise model above. 
 Hence, we characterize all systems explaining the data in a QMI of the form \eqref{e:Ineq_noise}, where $N \in \R^{14 \times 14}$.

First, we check the generalized Slater condition \eqref{e: slater on N} by verifying that $N$ has $7$ positive eigenvalues. Next, we can verify that conditions (i), (ii) and (iii) of Theorem~\ref{thm:inf_bal} are satisfied, for each $\sigma \in \{0.002,0.005,0.01,0.03,0.05\}$, meaning that the data are informative for generalized Lyapunov balancing. Conditions (ii) and (iii) of Theorem~\ref{thm:inf_bal} are semidefinite programs that we solve in Matlab, using Yalmip \cite{Yalmip} with SDPT3 \cite{toh1999sdpt3} as an LMI solver. 
As a result, we obtain generalized Hankel singular values that are common to all systems explaining the data. They are depicted in Figure~\ref{fig:Hankel} for various noise levels.

\begin{figure}[t!]
	\input{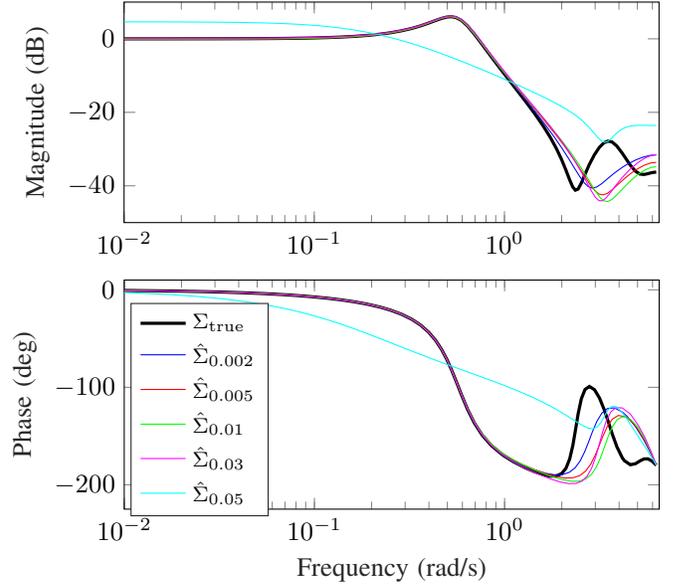}
	\caption{Bode plots of reduced-order models for five different noise levels. The systems $\hat{\Sigma}_{0.002}$, $\hat{\Sigma}_{0.005}$, $\hat{\Sigma}_{0.01}$, $\hat{\Sigma}_{0.03}$ and $\hat{\Sigma}_{0.05}$ denote reduced-order models of systems explaining the data with noise levels $\sigma=0.002,0.005,0.01,0.03$ and $0.05$, respectively, while  $\Sigma_{\mathrm{true}}$ denotes the true system. We stress that $\Sigma_{\sigma}$ is one from infinitely many systems in $\sigredbt_\sigma$ and $\Sigma_{\mathrm{true}}$ is assumed to be unknown.}
	\label{fig:bode}
\end{figure}

Figure~\ref{fig:Hankel} shows that the generalized Hankel singular values (generalized Gramians) obtained from the data indeed bound the ordinary Hankel singular values (Gramians) of the true system, see also the discussion in Section~\ref{section:GLB}. 
We stress however that the true system (and, hence, its Gramians) is assumed to be unknown. Additionally, we observe that the generalized Hankel singular values provide less strict bounds on the unknown ordinary Hankel singular values when the noise level is increased. 
%

In the balancing process, we also obtain the data-driven balancing transformation $T$ as well as the  data-driven projection matrices $\hatW$ and $\hatV$. Here, we take reduced models of order $r=3$ and therefore we have $\hatW,\hatV \in \R^{6 \times 3}$. From these projection matrices, the set $\sigredbt$, i.e., the set of reduced-order models via data-driven generalized balanced truncation, can be defined in terms of a QMI as stated in Corollary~\ref{cor:red_in_QMI}. 
We note that the set of reduced-order models $\sigredbt$ is only characterized by matrix $N_{V,W} \in \R^{8 \times 8}$, which is of reduced dimension (with respect to $N$). 

\begin{figure}[t!]
	\input{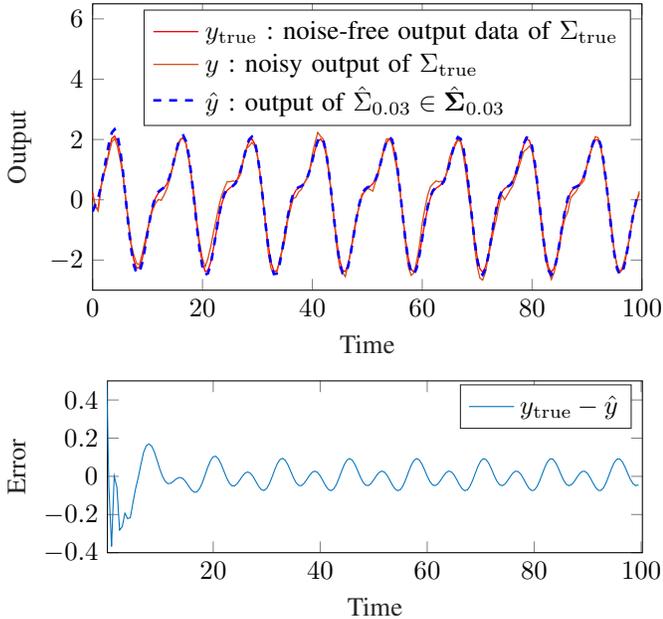}
%
%
\definecolor{mycolor1}{rgb}{0.00000,0.44700,0.74100}%
\begin{tikzpicture}

\begin{axis}[%
width=2.8in,
height=0.9in,
at={(0in,0in)},
scale only axis,
xmin=0.230414746543772,
xmax=100.230414746544,
xlabel style={font=\color{white!15!black}},
xlabel={Time},
ymin=-0.4,
ymax=0.5,
ylabel style={font=\color{white!15!black}},
ylabel={Error},
axis background/.style={fill=white},
legend style={legend cell align=left, align=left, draw=white!15!black}
]
\addplot [color=mycolor1]
  table[row sep=crcr]{%
0	0.621828123440507\\
0.5	-0.0557421781999294\\
1	-0.367152082222936\\
1.5	0.00217794721502196\\
2	-0.0572738928204369\\
2.5	-0.280956383448233\\
3	-0.263341512515807\\
3.5	-0.193024434452637\\
4	-0.221533448212171\\
4.5	-0.216864236000686\\
5	-0.138423258526087\\
5.5	-0.0640087239166971\\
6	-0.00500637014755356\\
6.5	0.0638645248337191\\
7	0.126589397275622\\
7.5	0.160721734810816\\
8	0.169278710726346\\
8.5	0.158456606800816\\
9	0.128386383809137\\
9.5	0.0850872158855258\\
10	0.0400373440199093\\
10.5	0.00173544291734229\\
11	-0.0251035177795474\\
11.5	-0.0377917694291151\\
12	-0.0370481267124272\\
12.5	-0.0274788636969313\\
13	-0.0154370419244105\\
13.5	-0.00694903031747596\\
14	-0.00649367937796563\\
14.5	-0.0158503122545344\\
15	-0.033423092535342\\
15.5	-0.0546898642006619\\
16	-0.0734833713675851\\
16.5	-0.0835897853060827\\
17	-0.0803892445676062\\
17.5	-0.0621856447060334\\
18	-0.0308130477003153\\
18.5	0.00865239186757406\\
19	0.0489758410764315\\
19.5	0.0823889003462595\\
20	0.10238211603308\\
20.5	0.10522432701672\\
21	0.0908319053297908\\
21.5	0.0627794794609899\\
22	0.0274500399141122\\
22.5	-0.00747343303968639\\
23	-0.0348106113792566\\
23.5	-0.0495436310100052\\
24	-0.0499340576476712\\
24.5	-0.0378307734909203\\
25	-0.0180998498152927\\
25.5	0.00267790839246451\\
26	0.0179169257424375\\
26.5	0.022720618885357\\
27	0.0151404236306787\\
27.5	-0.00331088425184611\\
28	-0.0280187291134486\\
28.5	-0.0524111732398118\\
29	-0.0696067480624185\\
29.5	-0.074159678782002\\
30	-0.0634715218139339\\
30.5	-0.0385189516480434\\
31	-0.00371417018769946\\
31.5	0.0340736277404795\\
32	0.0671019969545616\\
32.5	0.0885620599562342\\
33	0.0941711122533722\\
33.5	0.0831653468981584\\
34	0.0584528307332159\\
34.5	0.0258912037096775\\
35	-0.00713250933096377\\
35.5	-0.0334731379479452\\
36	-0.0479182483147015\\
36.5	-0.0483881888162697\\
37	-0.0363560482721512\\
37.5	-0.0163900586889614\\
38	0.00506369406207546\\
38.5	0.0213542603453685\\
39	0.0273391838337704\\
39.5	0.0207135845403357\\
40	0.0026253015369635\\
40.5	-0.0225776191073315\\
41	-0.0484070244285555\\
41.5	-0.0678430877007385\\
42	-0.0750960814931121\\
42.5	-0.0670875342993709\\
43	-0.0442842533979071\\
43.5	-0.0106761053427381\\
44	0.0271026502700765\\
44.5	0.0612993786293423\\
45	0.0848488674968184\\
45.5	0.0930224777303552\\
46	0.0845218634614322\\
46.5	0.0617534487813276\\
47	0.0302191603203701\\
47.5	-0.0028261505763747\\
48	-0.030122140943633\\
48.5	-0.0461291926290724\\
49	-0.0483139343519489\\
49.5	-0.0376836951635532\\
50	-0.0184443400935981\\
50.5	0.00312882514152146\\
51	0.0203244821400923\\
51.5	0.0277271228744994\\
52	0.0226177524302448\\
52.5	0.00569588130926335\\
53	-0.0190574654919626\\
53.5	-0.0453443220586693\\
54	-0.0661064174969965\\
54.5	-0.0752918566167127\\
55	-0.0693969162537251\\
55.5	-0.0484017495536229\\
56	-0.0158660754011443\\
56.5	0.0218442128728135\\
57	0.0570162619689807\\
57.5	0.0823842043287\\
58	0.0928244842683665\\
58.5	0.0865426890507184\\
59	0.0654655128399639\\
59.5	0.0347448367973169\\
60	0.00149685367857644\\
60.5	-0.0269150811679981\\
61	-0.0446350120418719\\
61.5	-0.0486797747406307\\
62	-0.0395851678173104\\
62.5	-0.0211799305222084\\
63	0.000451239986712648\\
63.5	0.0185555526363548\\
64	0.0274560026475733\\
64.5	0.0240186736882031\\
65	0.00847940291862792\\
65.5	-0.0155762811187756\\
66	-0.0420773809028698\\
66.5	-0.0639650712953213\\
67	-0.0749563878213817\\
67.5	-0.0711399478145962\\
68	-0.0520094469153596\\
68.5	-0.020677179811025\\
69	0.016790259887765\\
69.5	0.0527575310340718\\
70	0.079794319267902\\
70.5	0.0924089984249274\\
71	0.0883230374309223\\
71.5	0.0689784950004499\\
72	0.0391574824345968\\
72.5	0.00580847225941761\\
73	-0.023633007808889\\
73.5	-0.0430195052099497\\
74	-0.0489308022596111\\
74.5	-0.0414288489423549\\
75	-0.0239487614631492\\
75.5	-0.00235840809145482\\
76	0.0165769924987644\\
76.5	0.0269432008700574\\
77	0.0252044438693024\\
77.5	0.0111311352088126\\
78	-0.0121031665215905\\
78.5	-0.0386818349541502\\
79	-0.061577247137151\\
79.5	-0.0743037332936767\\
80	-0.0725602111711208\\
80.5	-0.055356928902849\\
81	-0.0253458688865782\\
81.5	0.0117315580366141\\
82	0.0483503349711674\\
82.5	0.0769476769603714\\
83	0.091687445640968\\
83.5	0.0898165733295593\\
84	0.0722859023697526\\
84.5	0.0434867449489778\\
85	0.0101695679406637\\
85.5	-0.0201900699384605\\
86	-0.0411794148929634\\
86.5	-0.0489551442639476\\
87	-0.043103791929085\\
87.5	-0.0266503192361911\\
88	-0.00521849420983983\\
88.5	0.0144451968153586\\
89	0.0262165397621038\\
89.5	0.0261744552777662\\
90	0.0136256191499631\\
90.5	-0.00868203892262631\\
91	-0.0352119502783568\\
91.5	-0.0589991074116718\\
92	-0.0733845354209586\\
92.5	-0.0736975604975765\\
93	-0.0584707693636319\\
93.5	-0.0298857549197546\\
94	0.00666486239415803\\
94.5	0.0437978833349351\\
95	0.0738498486668628\\
95.5	0.0906645348736337\\
96	0.0910255058103173\\
96.5	0.075387626599924\\
97	0.047731920353246\\
97.5	0.01458140391121\\
98	-0.0165805507960632\\
98.5	-0.0391030484266251\\
99	-0.0487351943627274\\
99.5	-0.044588309909103\\
};
\addlegendentry{$y_{\mathrm{true}}-\haty$}

\end{axis}

\begin{axis}[%
width=2.6in,
height=1in,
at={(0in,0in)},
scale only axis,
xmin=0,
xmax=1,
ymin=0,
ymax=1,
axis line style={draw=none},
ticks=none,
axis x line*=bottom,
axis y line*=left,
legend style={legend cell align=left, align=left, draw=white!15!black}
]
\end{axis}
\end{tikzpicture}%
	\caption{Time-domain output of a reduced-order model corresponding to the noise level $\sigma=0.03$ with input signal \eqref{e:input_signal}: $u(k)=2\sin(k)+\cos(0.5 k)$ compared to the noise-free and noisy output of the true system with the same input signal.}
	\label{fig:time_domain_output}
\end{figure}

In this example, for each noise level which is indicated by $\sigma$, we have a different set of reduced-order models denoted by $\sigredbt_\sigma$.
Then, from each set $\sigredbt_\sigma$, we pick a reduced-order model $\hat{\Sigma}_\sigma$. 
We stress that for a given noise level, $\hat{\Sigma}_\sigma$ depicts one from infinitely many possible reduced-order models contained in $\sigredbt_\sigma$. The Bode diagram of the reduced-order systems $\hat{\Sigma}_\sigma$'s is depicted in Figure~\ref{fig:bode}. Additionally, the time-domain output of  $\hat{\Sigma}_{0.03}$ is shown in Figure~\ref{fig:time_domain_output}. The Bode diagram in Figure~\ref{fig:bode} shows that reduced-order models accurately approximate the true system at least up to the noise level $\sigma=0.03$. But, if we increase the noise, e.g., $\sigma=0.05$, the resulting reduced-order model may not be able to accurately approximate the true system. From Figure~\ref{fig:time_domain_output}, we see that the reduced-order model corresponding to noise $\sigma=0.03$ is able to reconstruct the output data of the true system.

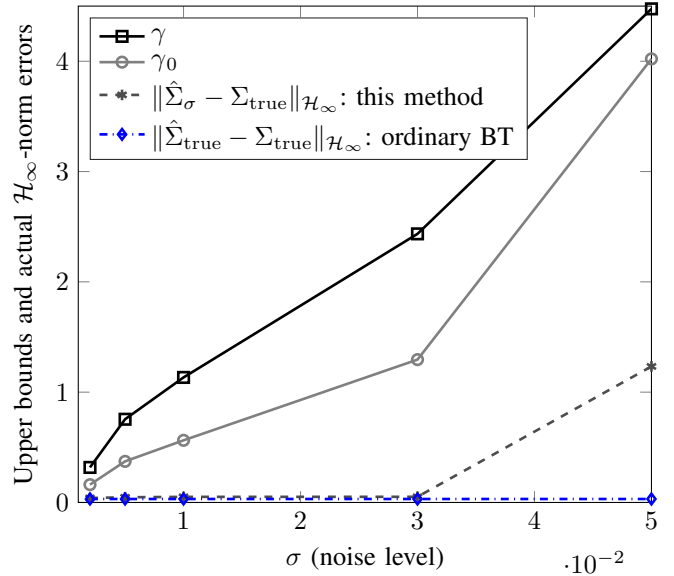
\begin{figure}[t!]
%
%
\begin{tikzpicture}

\begin{axis}[%
width=3in,
height=2.6in,
at={(0.011in,0.25in)},
scale only axis,
xmin=0.001,
xmax=0.05,
xlabel={$\sigma$ (noise level)},
ymin=0,
ymax=4.5,
ylabel={Upper bounds and actual $\HHH_\infty$-norm errors},
axis background/.style={fill=white},
legend style={at={(0.02,0.69)}, anchor=south west, legend cell align=left, align=left, draw=white!15!black}
]
\addplot [color=black, line width=1.0pt, mark=square, mark options={solid, black}]
  table[row sep=crcr]{%
0.002	0.31761\\
0.005	0.75431\\
0.01	1.13351\\
0.03	2.43501\\
0.05	4.47481\\
};
\addlegendentry{$\gamma$}

\addplot [color=gray, line width=1.0pt, mark=o, mark options={solid, gray}]
  table[row sep=crcr]{%
0.002	0.16151\\
0.005	0.37251\\
0.01	0.56321\\
0.03	1.29501\\
0.05	4.02141\\
};
\addlegendentry{$\gamma{}_{0}$}

\addplot [color=white!30!black, dashed, line width=1.0pt, mark=asterisk, mark options={solid, white!30!black}]
  table[row sep=crcr]{%
0.002	0.0405005870300648\\
0.005	0.0469910990332989\\
0.01	0.0507436403899286\\
0.03	0.0513144729349921\\
0.05	1.23313547638487\\
};
\addlegendentry{$\lVert \hat{\Sigma}_\sigma - \Sigma_{\mathrm{true}} \rVert_{\HHH_\infty}$: this method}

\addplot [color=blue, dashdotted, line width=1.0pt, mark=diamond, mark options={solid, blue}]
  table[row sep=crcr]{%
0.002	0.0314081111423348\\
0.005	0.0314081111423348\\
0.01	0.0314081111423348\\
0.03	0.0314081111423348\\
0.05	0.0314081111423348\\
};
\addlegendentry{$\lVert \hat{\Sigma}_{\mathrm{true}}- \Sigma_{\mathrm{true}} \rVert_{\HHH_\infty}$: ordinary BT}

\end{axis}

\begin{axis}[%
width=3.778in,
height=2.533in,
at={(0in,0in)},
scale only axis,
xmin=0,
xmax=1,
ymin=0,
ymax=1,
axis line style={draw=none},
ticks=none,
axis x line*=bottom,
axis y line*=left,
legend style={legend cell align=left, align=left, draw=white!15!black}
]
\end{axis}
\end{tikzpicture}%
	\caption{Comparison of  a priori error bounds $\gamma$ (solution of Problem~\eqref{prob:opt}), a posteriori error bounds $\gamma_0$ (solution of Problem~\ref{prob:opt0}), the actual $\HHH_\infty$-norm of the error between the true system $\Sigma_{\mathrm{true}}$ and a reduced-order model $\hat{\Sigma}_\sigma$ using this method and the $\HHH_\infty$-norm of the error between the tue system and its reduced-order model via ordinary balanced truncation (ordinary BT).}
	\label{fig:errorbound}
\end{figure}

Next, we will compute the error bounds for the reduced-order models in this framework. We note first that from the result of Theorem~\ref{thm:inf_bal}, all systems in $\sig$ and therefore $\sigredbt$ are guaranteed to be asymptotically stable. As a consequence, the LMIs \eqref{e:LMIs} and \eqref{e:LMI_0} are guaranteed to be feasible for some large enough $\gamma$ and $\gamma_0$, respectively. Hence, we can solve problems \eqref{prob:opt} and \eqref{prob:opt0}. The results can be found in 
Figure~\ref{fig:errorbound}. 

From 
Figure~\ref{fig:errorbound}, the upper bounds on the error with respect to the true system either from a priori or a posteriori upper bounds are getting more conservative when the noise levels are increased. This conservatism is caused  by the fact that the only knowledge that is available on the true system is that it is contained in $\sig$, the set of systems explaining the data. As the set $\sig$ has a larger size for increasing noise level, this leads to more conservative results.
 It is also clear that the a posteriori upper bound (the solution of problem \eqref{prob:opt0}) is less conservative than the corresponding a priori upper bound (solutions of problem~\eqref{prob:opt}) for each noise level.
 
In spite of the conservatism, the actual $\calH_\infty$-norms of the errors between the true system and the reduced-order models selected from $\sigredbt_\sigma$ for some small enough noise levels show that this data-driven method performs well.
In particular, the $\calH_\infty$-norm of the errors for noise levels $\sigma=0.002,0.005,0.01$ and $0.03$ which are given by $0.0 405$, $0.0470$, $0.0507$ and $0.0513$, respectively, are relatively small compared to the error of reduction by the ordinary balanced truncation, which is equal to $0.0314$. We stress, however, that computation via the ordinary balanced truncation requires the knowledge of the true system which cannot be achieved on the basis of the data.

\section{Conclusion} \label{sec:conclusion}
In this paper, a data-driven procedure to obtain reduced-order models from noisy data is developed. The procedure begins with introducing the concept of data reduction. Based on the noise model introduced by \cite{van2020noisy}, all (higher-order) systems explaining the data can be characterized in a high-dimensional quadratic matrix inequality (QMI) with a special structure. Due to this special structure, the class of reduced-order models obtained by applying a Petrov-Galerkin projection to all systems explaining the data can be characterized in a reduced-order QMI. 
As these QMIs depend only on the data and projection matrices, this can be regarded as a data reduction procedure.
Since this concept holds for general Petrov-Galerkin projections, it can potentially be extended to  solve data-driven model reduction problems via any projection-based reduction technique.

We then follow up the data reduction concept by constructing specific projection matrices from data. In particular, based on generalized controllability and observability Gramians, we provide necessary and sufficient conditions  such that all systems explaining the data have common generalized Gramians. 
These conditions 
substitute Lyapunov inequalities by a data-guided linear matrix inequality which can be solved efficiently by modern LMI solvers. Subsequently, a common balancing transformation and therefore common projection matrices  for generalized balanced truncation (which are in the class of the Petrov-Galerkin projections) are available to apply the data reduction. As such, a set of reduced-order models via generalized balanced truncation can then be characterized in a lower-dimensional QMI. Moreover, all reduced-order models in this set are guaranteed to be asymptotically stable and computable a priori and a posteriori upper bounds on the reduction error with respect to the true system are available. 

Beside the extension on exploiting data reduction via any projection-based reduction technique as mentioned above, ideas for future work include several directions. First, we aim at extending this result for input-output noisy data, even though an obstacle in this setting may be the construction of a state sequence. Second, investigating model reduction with preserving specific system properties such as network structure and port-Hamiltonian structure is often desirable. 



%

\appendix[Proofs]

\subsection{Proof of Proposition~\ref{l:N nonsingular}} \label{app:N_nonsingular}
\begin{proof}
	To prove the {`only if'} part, suppose that $\sig$ is bounded and has nonempty interior. Since $\sig$ has nonempty interior, then \eqref{e: slater on N} holds for some $\barS$. Next,
	let $N$ be defined as in \eqref{e: N} and partitioned as 
	$$ N=\bbma N_{11} & N_{12} \\ N_{12}^\top & N_{22} \ebma $$ 
	where $N_{11} \in \R^{(n+p) \times (n+p)}$, $N_{12} \in \R^{(n+p) \times (n+m)}$, and $N_{22} \in \R^{(n+m) \times (n+m)}$. 
	By Lemma~\ref{l:bounded_nonempty}, matrix $N$ is nonsingular with ${N_{11}-N_{12}N_{22}^{-1}N_{12}^\top >0}$ and $N_{22}<0$. As a result, the matrix 
	$$ \begin{bmatrix}
	I & 0  & X_{+} \\0 & I & Y _{-}\\ 0 &0 & -X_{-} \\ 0 & 0 & -U_{-}
	\end{bmatrix}$$
	and therefore $\bbma X_{-} \\ U_{-} \ebma$ are full row rank.
	
	To prove the {`if'} part, observe that it follows from \eqref{e: slater on N} and \cite[Fact~5.8.16]{bernstein2009matrix} that
	\beq \label{e:upsilon_+}\upsilon_+(N) \geq n+p. \eeq
	Note that 
	$$ N_{22}= \bbma X_{-} \\ U_{-} \ebma \Phi_{22} \bbma X_{-} \\ U_{-} \ebma^\top. $$
	Since $\Phi_{22} <0$ and  $\bbma X_{-}  \\  U_{-} \ebma$ has full row rank, we have that $N_{22} <0$. Then, following Haynsworth's inertia additivity formula, see \cite[Fact 6.5.5]{bernstein2009matrix}, we have 
	$$ \In(N)=\In(N_{22}) + \In(N_{11}-N_{12}N_{22}^{-1}N_{12}^\top),$$
	and hence
	\beq \label{e:upsilon_-} \upsilon_-(N) \geq n+m. \eeq
	Since $N \in \R^{(m+2n+p) \times (m+2n+p)}$, $\upsilon_-(N) + \upsilon_+(N) \leq m+2n+p$. This, together with \eqref{e:upsilon_+} and \eqref{e:upsilon_-}, implies that $\upsilon_+(N)+\upsilon_-(N)=m+2n+p$. Therefore, $N$ has no zero eigenvalues. Consequently, it is nonsingular with $N_{22}<0$ and $N_{11}-N_{12}N_{22}^{-1}N_{12}^\top >0$. Then, it follows from \eqref{e:Ineq_noise} and Lemma~\ref{l:bounded_nonempty} that $\sig$ is bounded and has nonempty interior. 
\end{proof}
\subsection{Proof of Theorem~\ref{thm:ROM in QMI}} \label{app_proj}
Before giving the proof of Theorem~\ref{thm:ROM in QMI}, we will develop some general results on quadratic matrix inequalities. 

To this end, consider the set 
\beq \label{setN}
\calM:=\set{Z \in\R^{p\times q}}{\bbm I\\Z^\top \ebm^\top \Psi \bbm I\\Z^\top \ebm\geq 0},
\eeq
where $\Psi$ admits the partitioning
\beq \label{e:partitionPsi}
\Psi= \bbm \Psi_{11}&\Psi_{12}\\\Psi_{12}^\top &\Psi_{22}\ebm
\eeq 
with $\Psi_{11} \in \R^{p \times p}$. 
We assume throughout this appendix that $\Psi_{22} < 0$ and $\Psi_{11}-\Psi_{12}\Psi_{22}^{-1}\Psi_{12}^\top >0$. It is known from Lemma~\ref{l:bounded_nonempty} that $\calM$ is bounded and has nonempty interior. As a consequence of this, $\calM$ admits various representations, as stated next.
\blemp \label{l:setN_equiv}
Consider $\calM$ given by \eqref{setN}.
Define $\Psi | \Psi_{22} := \Psi_{11}-\Psi_{12}\Psi_{22}^{-1}\Psi_{12}^\top$ and 
\beq \label{e:Xi}
\Psi^{\sharp} :=\bbm 0 & -I_q \\I_p & 0\ebm
\Psi\inv
\bbm 0 & -I_p\\I_q & 0\ebm.
\eeq 
Suppose that $\Psi | \Psi_{22} > 0$ and $\Psi_{22}<0$. Then, $\calM=\calM_1=\calM_2$, where
\begin{align}
	 &\calM_1 :=\set{Z}{\bbm
		I \\ Z
		\ebm^\top
		\Psi^{\sharp}
		\bbm
		I \\ Z
		\ebm\geq 0 }
	, \nonumber \\
	&\calM_2 := \set{ Z } {\bbm I \\ Z^\top \ebm^\top \!\! \Psi \!\bbm I \\ Z^\top \ebm \!=\! \calQ \text{ where } 0\leq \calQ \leq \Psi | \Psi_{22}\!}\!. \nonumber
\end{align}
\elemp
\begin{proof}
	The proof of  $\calM =   \calM_1$ is provided in \cite{vanDissipativity2022mag}. We will prove that $\calM =   \calM_2$.
	Clearly, $ \calM_2\subseteq \calM$. Then, it remains to show that the reverse inclusion holds. Let $Z \in \calM$  and let 
	$$
	\calQ = \bbm I \\ Z^\top \ebm^\top \Psi \bbm I \\ Z^\top \ebm.
	$$
	Clearly, $\calQ \geq 0$. Note that
	$$
	\calQ =\Psi|\Psi_{22}+ (Z^\top+\Psi_{22}\Psi_{12}^\top)^\top \Psi_{22} (Z^\top +\Psi_{22}\Psi_{12}^\top) \leq \Psi|\Psi_{22}
	$$
	since $\Psi_{22} <0$. Therefore, $0\leq \calQ \leq \Psi |\Psi_{22}$ and, as a result, ${Z \in \calM_2}$.
\end{proof}
\begin{figure*}[h]
	\centering
	\normalsize
	\begin{equation}\label{e:PsiWV} 
	\tag{28} 
	\nonumber
	\Psi_{V,W} :=  \begin{pmat}[{|}]
	W^\top (\Psi | \Psi_{22} +\Psi_{12}\Psi_{22}^{-1} V (V^\top \Psi_{22}^{-1}V)^{-1}V^\top \Psi_{22}^{-1}\Psi_{12}^\top )W &  W^\top \Psi_{12} \Psi_{22}^{-1} V (V^\top \Psi_{22}^{-1}V)^{-1} \cr\-
	(V^\top \Psi_{22}^{-1}V)^{-1}V^\top \Psi_{22}^{-1}\Psi_{12}^\top W  & (V^\top \Psi_{22}^{-1} V)^{-1} \cr
	\end{pmat}
	\end{equation}
	\begin{equation}\label{e:PsiW} \tag{29} 
	\nonumber
	{\Psi}_{V} :=  \begin{pmat}[{|}]
	\Psi | \Psi_{22} +\Psi_{12}\Psi_{22}^{-1} V (V^\top \Psi_{22}^{-1}V)^{-1}V^\top \Psi_{22}^{-1}\Psi_{12}^\top  &   \Psi_{12} \Psi_{22}^{-1} V (V^\top \Psi_{22}^{-1}V)^{-1} \cr\-
	(V^\top \Psi_{22}^{-1}V)^{-1}V^\top \Psi_{22}^{-1}\Psi_{12}^\top   & (V^\top \Psi_{22}^{-1} V)^{-1} \cr
	\end{pmat}
	\end{equation}
	\begin{equation}\label{e:XiW} \tag{30} 
	\nonumber
	{\Psi}_{V}^{\sharp} :=  \begin{pmat}[{|}]
	- V^\top (\Psi_{22}^{-1} - \Psi_{22}^{-1}\Psi_{12}^\top (\Psi | \Psi_{22})^{-1} \Psi_{12} \Psi_{22}^{-1})V  &   -V^\top \Psi_{22}^{-1}\Psi_{12}^\top (\Psi | \Psi_{22})^{-1} \cr\-
	-(\Psi | \Psi_{22})^{-1} \Psi_{12} \Psi_{22}^{-1}V   &  -(\Psi | \Psi_{22})^{-1}\cr 
	\end{pmat}
	\end{equation}
	\par\noindent\rule{\textwidth}{0.5pt}
\end{figure*}
Now, we will consider projections of the elements of the set $\calM$. Let $V \in \R^{q \times \hat{q}}$ and ${W \in \R^{p \times \hat{p}}}$  be full column rank projection matrices with $\hatp \leq p$ and $\hatq \leq q$. We refer to the set
$$
\calM_{V,W}:= \set{W^\top Z V }{Z \in \calM}
$$
as a reduction of  $\calM$ using projectors $W$ and $V$. Note that we do not assume that $W^\top V =I$.

We will show that elements of $\calM_{V,W}$ themselves satisfy a quadratic matrix inequality. To do so, let $Z \in \calM$ such that we have
\beq \nonumber
\bbm I\\Z^\top \ebm^\top \bbm \Psi_{11}&\Psi_{12}\\\Psi_{12}^\top &\Psi_{22}\ebm\bbm I\\Z^\top \ebm\geq 0,
\eeq  
which can be written as
\newcounter{counterafterM}
\normalsize
\setcounter{counterafterM}{\value{equation}}
\setcounter{equation}{30}
\beq 
\label{e:QMI_setN}
\Psi | \Psi_{22} + (Z + \Psi_{12} \Psi_{22}^{-1})\Psi_{22}(Z^\top+ \Psi_{22}^{-1}\Psi_{12}^\top) \geq 0.
\eeq
Using the Schur complement, \eqref{e:QMI_setN} is equivalent to
\beq \label{e:Schur_of_QMIsetN}
\bbm \Psi | \Psi_{22} & Z + \Psi_{12} \Psi_{22}^{-1} \\
Z^\top+ \Psi_{22}^{-1}\Psi_{12}^\top & -\Psi_{22}^{-1} 
\ebm \geq 0,
\eeq 
where we note that the inverse of $\Psi_{22}$ exists as $\Psi_{22}<0$ by assumption. Pre- and postmultiplying \eqref{e:Schur_of_QMIsetN} by $\blkdiag(W^\top,V^\top)$  and $\blkdiag(W,V)$, respectively, gives
\beq \label{e:projSchur_setNr}
\bbm W^\top (\Psi | \Psi_{22})W& W^\top (Z + \Psi_{12} \Psi_{22}^{-1}) V \\
V^\top (Z^\top+ \Psi_{22}^{-1}\Psi_{12}^\top)W & -V^\top \Psi_{22}^{-1} V
\ebm \geq 0.
\eeq  
Let us define $\hat{Z}=W^\top Z V$. Using again a Schur complement argument and writing the result as a quadratic matrix inequality, we obtain that \eqref{e:projSchur_setNr} is equivalent to
\beq \label{e:QMI_setNr}
\bbm I\\ \hatZ^\top \ebm^\top  \Psi_{V,W}  \bbm I\\ \hatZ^\top \ebm\geq 0,
\eeq
 where $\Psi _{V,W}$ is given by \eqref{e:PsiWV}. 
This shows that if $Z \in \calM$ then $W^\top Z V$ satisfies \eqref{e:QMI_setNr}. 
Stated differently, we have that $\calM_{V,W} \subseteq \set{\hatZ}{\eqref{e:QMI_setNr} \text{ holds}}$. In fact, we have that the equality holds as asserted in the following theorem. 
\bthep\label{thm:Zr_eq}
It holds that
\beq  \nonumber
\begin{split}
	\calM_{V,W}&= \set{\hatZ}{  
		\bbm I\\ \hatZ^\top \ebm^\top  \Psi_{V,W}  \bbm I\\ \hatZ^\top \ebm\geq 0 
	}.
\end{split}
\eeq
\ethep
Before we present the proof of Theorem~\ref{thm:Zr_eq}, we need to state some auxiliary lemmas.
First, we recall the following result from linear algebra \cite[Fact~5.10.19]{bernstein2009matrix}.
\blemp \label{l:XX=YY}
	Let $X,Y \in \R^{p \times q} $. Then, $X ^\top X = Y^\top Y$ if and only if  $Y = U X$ where $U$ is an orthogonal matrix, i.e., ${U^\top U =I}$.
\elemp
 The next two lemmas are central to prove Theorem~\ref{thm:Zr_eq}.
\blemp \label{l:Q_full}
Let $S=S^\top \in \R^{p \times p}$ such that $S>0$. Let $V \in \R^{p \times r}$ be a full rank matrix with $r \leq  p$ and $Q_V=Q_V^\top \in \R^{r \times r}$ such that $0 \leq  {Q}_V  \leq  V^\top S V$. Then, there exists $Q =Q^\top$ such that 
$
0 \leq Q  \leq S  
$ and $V^\top Q V = Q_V$.
\elemp 
\begin{proof}
	Since $V$ is full rank, there exists a matrix $\tilde{V}$ such that $T_V := \bbm V & \tilde{V} \ebm$ is nonsingular. Note that  
	\beq \nonumber
	T_V^\top  S T_V =F^\top  \bbm V^\top S  V & 0 \\ 0 &  (T_V^\top S T_V) | (V^\top S  V) \ebm F
	\eeq
	where 
	$$
	F := \bbm I &  (V^\top S  V)^{-1}V^\top S \tilde{V}  \\ 0& I \ebm .
	$$
	Let 
	\beq \nonumber 
	\bar{Q} :=  F^\top \bbm Q_V & 0 \\ 0 & \Delta \ebm F,
	\eeq
	where $0\leq \Delta \leq (T_V^\top S T_V) | (V^\top S  V)$. Clearly, we have ${0 \leq \bar{Q}  \leq T_V^\top S T_V}$. Take $Q = T_V^{-\top} \bar{Q} T_V^{-1}$ to guarantee that $0 \leq Q \leq S$ and $V^\top Q V = Q_V$.
\end{proof}
\blemp \label{l:ZV}
Consider the matrix $\Psi \in \R^{(p+q) \times (p+q)}$ partitioned as in \eqref{e:partitionPsi} such that $\Psi_{11}-\Psi_{12}\Psi_{22}^{-1}\Psi_{12}^\top >0$ and $\Psi_{22}<0$. Let $W \in \R^{p \times r}$ be a full rank matrix with $r \leq p$. Suppose that $Z_W$ satisfies
\beq \label{e:Vphi}
\bbm I \\ Z_W^\top \ebm^\top  
\bbm W^\top \Psi_{11}  W & W^\top \Psi_{12} \\ \Psi_{12}^\top W  & \Psi_{22} \ebm \bbm I \\ Z_W^\top \ebm   \geq 0.
\eeq 
Then, there exists $Z$ such that  $Z^\top W=Z_W^\top$ and
\beq \nonumber 
\bbm I \\ Z^\top \ebm^\top  
\Psi \bbm I \\ Z^\top \ebm   \geq 0.
\eeq 
\begin{proof}
	Note that \eqref{e:Vphi} can be written as
	$$
	W^\top (\Psi | \Psi_{22}) W + \barZ_W \Psi_{22} \barZ_W^\top \geq 0,
	$$
	where $\bar{Z}_W =Z_W+ W^\top \Psi_{12}\Psi_{22}^{-1}$.
	Let \beq \label{e:Qv} Q_W := W^\top (\Psi | \Psi_{22}) W + \barZ_W \Psi_{22} \barZ_W^\top .\eeq Then, $0 \leq Q_W \leq W^\top (\Psi | \Psi_{22}) W$ since $\Psi_{22}<0$. By Lemma~\ref{l:Q_full}, there exists $Q=Q^\top$ such that $0 \leq Q \leq \Psi | \Psi_{22}$ and $W^\top Q W = Q_W$. As $\Psi|\Psi_{22}-Q\geq 0$, there exists a matrix $R$ such that
	$$
	\Psi | \Psi_{22}-Q = R^\top R,
	$$
	which implies
	\beq \label{e:Qv_prime}
	Q_W=W^\top (\Psi | \Psi_{22}) W - W^\top R^\top R W.
	\eeq 
	By comparing \eqref{e:Qv} and \eqref{e:Qv_prime}, it follows from Lemma~\ref{l:XX=YY} that
	$$
	(-\Psi_{22})^{\frac{1}{2}} \barZ_W^\top= U RW
	$$
	for some orthogonal matrix $U$. Note that since $W$ is full rank, there exists $\tilde{W}$ such that $\bbm W & \tilde{W} \ebm$ is nonsingular. Now, let $\tilde{Z}_W$ be such that 
	$$
	(-\Psi_{22}) ^{\frac{1}{2}} \tilde{Z}_W^\top = U R\tilde{W}.
	$$
	Next, define $$Z^\top = \bbm \barZ_W^\top  & \tilde{Z}_W^\top \ebm\bbm W & \tilde{W} \ebm^{-1} - \Psi_{22}^{-1}\Psi_{12}^\top, $$
	which can be easily checked to verify $Z^\top W = Z_W^\top$.  
Moreover,
	\beq \nonumber
	\begin{split}
		&\bbm I \\ Z^\top \ebm^\top \bbm \Psi_{11}&\Psi_{12}\\\Psi_{12}^\top &\Psi_{22}\ebm \bbm I \\ Z^\top \ebm \\ &= \Psi | \Psi_{22} + (Z+\Psi_{12}\Psi_{22}^{-1}) \Psi_{22} (Z+\Psi_{12}\Psi_{22}^{-1})^\top \\ 
		&=\Psi | \Psi_{22}- R^\top R= Q \geq 0,
	\end{split}
	\eeq 
	as desired.
\end{proof}
\elemp 
Now, we are ready to prove Theorem~\ref{thm:Zr_eq}.
\begin{proof}[Proof of Theorem~\ref{thm:Zr_eq}]
	It is clear that $${\calM_{V,W} \subseteq \set{\hatZ}{\eqref{e:QMI_setNr}\text{ holds}}}.$$
	To prove the reverse inclusion, let $\hatZ$ be such that \eqref{e:QMI_setNr} holds, i.e.,
	$$
	\bbm I\\ \hatZ^\top \ebm^\top  \Psi_{V,W}  \bbm I\\\hatZ^\top \ebm\geq 0.
	$$
	We will show that there exists $Z \in \calM$ such that $\hatZ=W^\top Z V$. By Lemma~\ref{l:ZV}, there exists $Z_V$ such that $Z_V^\top W=\hatZ^\top $ and 
	\beq \label{e:QMI_XiV}
	\bbm I \\ Z_V^\top \ebm^\top {\Psi}_{V} \bbm I \\ Z_V^\top \ebm \geq 0,
	\eeq 
	where ${\Psi}_{V}$ is given by \eqref{e:PsiW}. From Lemma~\ref{l:setN_equiv}, \eqref{e:QMI_XiV} is equivalent to
	\beq \label{e:QMI_XiW}
	\bbm I \\ Z_V \ebm^\top {\Psi}^{\sharp}_{V} \bbm I \\ Z_V \ebm \geq 0,
	\eeq 
	where $\Psi_{V}^{\sharp}$ is given by \eqref{e:XiW}. Using Lemma~\ref{l:ZV} again, \eqref{e:QMI_XiW} implies the existence of $Z$ such that $Z V= Z_V$ and 
	\beq \nonumber \label{e:QMI_XiZ}
	\bbm I \\ Z \ebm^\top \Psi^{\sharp} \bbm I \\ Z \ebm \geq 0,
	\eeq 
	where $\Psi^{\sharp}$ is given by \eqref{e:Xi}. Therefore, there exists  $Z$ such that $\hatZ=W^\top Z V$ and, due to Lemma~\ref{l:setN_equiv}, $Z \in \calM$.
\end{proof}

\begin{proof}[Proof of  Theorem~\ref{thm:ROM in QMI}]
	Since \eqref{e: slater on N} holds for some $\barS$ and matrix $\bbma X_{-} \\ U_{-} \ebma$ has full row rank, $\sig$ is bounded and has nonempty interior. In addition, since ${\hatW^\top \hatV =I}$, then $W$ and $V$ are full column rank. Therefore, the claim follows from the result of Theorem~\ref{thm:Zr_eq}.
\end{proof}
\subsection{Strict matrix S-Lemma} \label{app:Slemma}
\newcounter{cprop} \counterwithin{cprop}{section}
\newtheorem{propositionp}[cprop]{Proposition}
\begin{propositionp}[{\cite[Thm. 11]{van2020noisy}}] \label{matSlemma_strict}
	Let $F,G \in \mathbb{R}^{(q+r) \times (q+r)}$ be symmetric matrices. Assume that 
	$$ 
	\Sigma_G : = \set{ {V} \in \mathbb{R}^{r \times q} } { \begin{bmatrix} I \\ V \end{bmatrix}^\top\!\!\! G \begin{bmatrix} I \\ V \end{bmatrix} \geq 0  }
	$$
	is bounded. Consider the statements
	\begin{enumerate}
		\renewcommand{\theenumi}{(\roman{enumi})}
		\renewcommand{\labelenumi}{\theenumi}
		\item\label{i:slemma1_s} There exists some matrix $\bar{V} \in \mathbb{R}^{r \times q}$ such that 
		\begin{equation} \nonumber 
		\label{matSlater}
		\begin{bmatrix} I \\ \bar{V} \end{bmatrix}^\top\!\!\! G \begin{bmatrix} I \\ \bar{V} \end{bmatrix} > 0.
		\end{equation}
		\item\label{i:slemma2_s} 
		$
		\begin{bmatrix} I \\ V \end{bmatrix}^\top\!\!\! F \begin{bmatrix} I \\ V \end{bmatrix} > 0 \:\: \forall\,  V \in \mathbb{R}^{r \times q} \text{ with} \begin{bmatrix} I \\ V \end{bmatrix}^\top\!\!\! G \begin{bmatrix} I \\ V \end{bmatrix}  \geq 0.
		$\\
		\item\label{i:slemma3_s} There exists a scalar $\alpha \geq 0$ such that $F - \alpha G > 0$.
	\end{enumerate}
	Then, the following implications hold:
	\ben
	\renewcommand{\theenumi}{(\Roman{enumi})}
	\renewcommand{\labelenumi}{\theenumi}
	\item\label{I:slemma1_s} \ref{i:slemma1_s} and \ref{i:slemma2_s} $\implies$ \ref{i:slemma3_s}.
	\item\label{I:slemma2_s} \ref{i:slemma3_s} $\implies$ \ref{i:slemma2_s}.
	\een
\end{propositionp}

\subsection{Proof of Theorem~\ref{thm:inf_bal}} \label{app:thm_inf}
\begin{proof}
	Let us first prove the `only if' statement. Suppose that the data $\dataUXY$  are informative for GLB. By Definition~\ref{def:Lyap_balancing}, there exist $P=P^\top>0$ and $Q=Q^\top>0$ such that \eqref{e:Ineq_ctrbGramian} and \eqref{e:Ineq_obsvGramian} hold for all $(A,B,C,D)$ satisfying \eqref{e:Ineq_noise}. We begin with statement \ref{mainthm1}. Let $\xi \in \R^n$ and $\eta \in \R^{m}$ be such that 
	\begin{equation*}
	\begin{bmatrix}
	\xi^{\top} & \eta^{\top}
	\end{bmatrix}\bbma X_{-} \\ U_{-} \ebma =0.
	\end{equation*}
	Moreover, let $(A,B,C,D) \in \sig$ and $\zeta \in \R^{n}$  be a nonzero vector. Note that
	\begin{equation*}
	(A+\alpha \zeta\xi^\top,B+\alpha \zeta \eta^\top,C,D) \in \sig
	\end{equation*}
	for every $\alpha \in \R$, as can be concluded from \eqref{e:Ineq_noise}. Since the data are informative for GLB, there exists $P=P^\top>0$ such that 
	\begin{equation}\label{e:lyap_Gramian_alpha}
	P-A_{\alpha}PA_{\alpha}^\top - B_{\alpha}B_{\alpha}^\top > 0
	\end{equation}
	where $A_{\alpha}:=A+\alpha \zeta \xi^\top$ and $B_{\alpha}=B+\alpha\zeta \eta^\top$. 
	Note that \eqref{e:lyap_Gramian_alpha} holds for every $\alpha \in \R$. Then,
	by dividing \eqref{e:lyap_Gramian_alpha} by $\alpha^2$ and letting $\alpha \to \infty$, we obtain
	$$
	(- \xi^\top P \xi  -  \eta^\top \eta) \zeta \zeta^\top \geq 0.
	$$
	Since $P>0$ and $\zeta\neq0$, we see that $\xi =0$ and $\eta =0$. Therefore, $\bbma X_{-} \\ U_{-} \ebma $ has full row rank.
	
	To show (ii) and (iii), we first rewrite the matrix inequalities \eqref{e:Ineq_ctrbGramian} and \eqref{e:Ineq_obsvGramian} as the quadratic matrix inequalities
	\beq \label{e:QMI_ctrb}
	\bbm I \\ A^\top \\ B^\top \ebm^\top   \bbm P & 0 & 0 \\ 0 & - P & 0 \\ 0& 0 & -I_m \ebm  \bbm I \\ A^\top \\ B^\top \ebm >0
	\eeq 
	and 
	\beq \label{e:QMI_obsv}
	\bbm I \\ A \\ C \ebm^\top   \bbm Q & 0 & 0 \\ 0 & -Q & 0 \\ 0& 0 & -I_p \ebm  \bbm I \\ A \\ C \ebm >0,
	\eeq 
respectively.
Note that  we have the QMI \eqref{e:Ineq_noise} characterizing the set of all systems explaining  the data. Moreover, since $\bbma X_{-} \\  U_{-} \ebma $ has full row rank and we assume that there exists $\barS$ such that \eqref{e: slater on N} holds, it follows from the proof of Proposition~\ref{l:N nonsingular} that systems explaining the data are equivalently characterized by
	\beq \label{e:QMI_noise_T}
	\systrr^\top \! N^{\sharp} \systrr \geq 0,
	\eeq
	where 
	\beq \nonumber 
	N^{\sharp} :=\begin{bmatrix} 0 & -I_{n+m} \\ I_{n+p} & 0 \end{bmatrix} N^{-1} \begin{bmatrix} 0 & -I_{n+p} \\ I_{n+m} & 0 \end{bmatrix},
	\eeq
	see Lemma~\ref{l:setN_equiv}. From a projection of  \eqref{e:Ineq_noise} and Lemma~\ref{l:ZV}, we have that
	all $(A,B)$ satisfying \eqref{e:Ineq_noise} are equivalent to those satisfying
	\beq \label{e:QMI_Nc}
	\bbm I \\ A^\top \\ B^\top  \ebm^\top N_{\calC}  \bbm I \\ A^\top \\ B^\top  \ebm \geq 0,
	\eeq 
	where $N_{\calC}$ is given by 
	\beq \nonumber
	N_{\calC}:= \bbm I_n & 0 \\ 0 & 0 \\ 0 & I_{n+m} \ebm^\top N \bbm I_n & 0 \\ 0 & 0 \\ 0 & I_{n+m} \ebm.
	\eeq 
	Similarly, all $(A,C)$ satisfying \eqref{e:QMI_noise_T} are equivalent to those satisfying
		\beq \label{e:QMI_No}
	\bbm I \\ A \\ C  \ebm^\top N_{\calO}  \bbm I \\ A \\ C  \ebm \geq 0,
	\eeq 
	where
	\beq \nonumber 
	{N}_{\calO} :=\bbm I_n & 0 \\ 0 & 0 \\ 0 & I_{n+p} \ebm^\top N^{\sharp} \bbm I_n & 0 \\ 0 & 0 \\ 0 & I_{n+p} \ebm.
	\eeq
	Now, we are ready to apply the matrix S-lemma from Appendix~\ref{app:Slemma}. In particular,  by informativity for GLB, \eqref{e:QMI_ctrb} holds for all $(A,B)$ satisfying \eqref{e:QMI_Nc}, such that the use of the matrix {S-lemma} (Proposition~\ref{matSlemma_strict} in Appendix~\ref{app:Slemma}) yields \eqref{e:Slemma_ctrb} and proves (ii). The proof of (iii) is similar, using \eqref{e:QMI_obsv} and \eqref{e:QMI_No}.
	
	
	To prove the `if' statement, first suppose that $\bbma X_- \\ U_- \ebma$ has full row rank. Then, under assumption \eqref{e: slater on N}, $N$ is nonsingular with $N_{22}<0$ and $N_{11}-N_{12}N_{22}^{-1}N_{12}^\top >0$. Thus, $N_{\calO}$ is well-defined. Now, suppose that statements \ref{mainthm2} and \ref{mainthm3} are satisfied. Then, the matrix S-lemma in Proposition~\ref{matSlemma_strict} implies that \eqref{e:QMI_ctrb} and \eqref{e:QMI_obsv} hold for all $(A,B)$ and $(A,C)$ satisfying \eqref{e:QMI_Nc} and \eqref{e:QMI_No}, respectively.  This implies that \eqref{e:Ineq_ctrbGramian} and \eqref{e:Ineq_obsvGramian} hold for all systems explaining the data, i.e., the data are informative for generalized Lyapunov balancing.
\end{proof} 

\subsection{Proof of Theorem~\ref{thm:maxG1G2}} \label{app:proofbound}
\begin{proof}
	We will prove the upper bound  by employing the bounded real lemma. 
	To do so, consider any $\Sigma \in \sig$ and $\hat{\Sigma} \in \sigredbt$ with realizations $(A,B,C,D)$ and $(\hatA,\hatB,\hatC,\hatD)$, respectively. Then,  a realization for $\hat{\Sigma}-\Sigma$ is given by the quadruplet
	\beq \nonumber
	A_\mathrm{d} \!:= \!\bbm A & 0 \\ 0 & \hatA \ebm , \ B_\mathrm{d} \!:= \!\bbm B \\ \hatB \ebm, \ C_\mathrm{d} := \!\bbm C & -\hatC \ebm, \ D_\mathrm{d} := \!D-\hatD.
	\eeq
	Let $\gamma >0$. By (the discrete-time version of) the bounded real lemma, e.g.,\cite[Thm. 4.6.6 (iv)]{skelton1997unified}, the matrix $A_{\mathrm{d}}$ satisifes $\rho(A_{\mathrm{d}}) <1$ and  ${\lVert \hat{\Sigma}- \Sigma \rVert_{\HHH_\infty} < \gamma}$ if and only if there exists ${K \in \R^{(n+r) \times (n+r)}}$ with $K=K^\top  >0$ such that
	\beq  \label{e:br_lemma}
	\bbm K & 0 \\ 0 & I_p \ebm - 
	\bbm A_{\mathrm{d}} & B_{\mathrm{d}} \\ C_{\mathrm{d}} & D_{\mathrm{d}} \ebm \bbm K & 0 \\ 0 & \gamma^{-2}I_m \ebm \bbm A_{\mathrm{d}} & B_{\mathrm{d}} \\ C_{\mathrm{d}} & D_{\mathrm{d}} \ebm^\top >0.
	\eeq
	If  \eqref{e:br_lemma} holds for all $\Sigma\in \sig$ and $\hat{\Sigma} \in \sigredbt$ (for the same $K$), then clearly the norm $\lVert \Sigma - \hat{\Sigma} \rVert_{\HHH_{\infty}}$  is upper bounded by $\gamma$ for all choices of systems in $\sig$ and $\sigredbt$.  
	Note that \eqref{e:br_lemma} can be written in the QMI form 
	\beq \label{e:br_in_QMI}
	 \bbm I  & 0\\ 0 & I \\ A_{\mathrm{d}}^\top & C_{\mathrm{d}}^\top \\ B_{\mathrm{d}}^\top & D_{\mathrm{d}}^\top \ebm^\top  
	\bbm K & 0 & 0&  0\\ 0 &  I_p & 0 & 0\\ 0 & 0 & -K & 0\\ 0 & 0& 0 & -\gamma^{-2} I_m \ebm 
	\bbm I  & 0\\ 0 & I \\ A_{\mathrm{d}}^\top & C_{\mathrm{d}}^\top \\ B_{\mathrm{d}}^\top & D_{\mathrm{d}}^\top \ebm \!> \!0.\!
	\eeq
	
	We will show the equivalence of \eqref{e:LMIs} and the satisfaction of \eqref{e:br_in_QMI} for all systems in $\sig$ and $\sigredbt$ by using the matrix S-lemma. As a first step, we introduce the notation
	\beq \label{e:Jnotation}
	J := \systr \text{ and } \hat{J} := \systrred,
	\eeq 
	such that the data equations \eqref{e:Ineq_noise} and
	\beq \nonumber 
   \systrred^\top \! \! \! \! N_{V,W} \systrred\!\geq 0
	\eeq 
	can be written as
	$
	J^\top {N} J \geq 0
	$ and 
	$
	\hat{J}^\top {N}_{V,W} \hat{J} \geq 0,
	$ respectively.
	On the other hand, it can be checked that  \eqref{e:br_in_QMI} is equivalent to 
	\beq \label{QMI_br_Gam}
	\Gamma^\top \bbm J &  0\\ 0 & \hat{J} \ebm ^\top \bbm \tilde{\Theta}_{11} & \tilde{\Theta}_{12} \\ \tilde{\Theta}_{12}^\top & \tilde{\Theta}_{22} \ebm \bbm J & 0 \\ 0 & \hat{J} \ebm \Gamma >0,
	\eeq
	where 
	\beq \nonumber 
	\tilde{\Theta}_{12}:= 
	\blkdiag(K_{12},0,-K_{12},-\gamma^{-2}I_m),
	\eeq
	\beq \nonumber
	\tilde{\Theta}_{ii}:= 
	\blkdiag(K_{ii},\frac{1}{2}I_p,-K_{ii},-\gamma^{-2}I_m)
	\eeq 
	for $ i=1,2$, and 
	$$
	\Gamma := \bbm I_n & 0 & 0 \\0 & 0 & I_p \\0 &I_r &0 \\0 & 0 &- I _p
	\ebm.
	$$
	Note that \eqref{QMI_br_Gam} means that
	\beq \nonumber 
	x^\top \bbm J &0 \\  0 & \hat{J} \ebm^\top  \bbm \tilde{\Theta}_{11} & \tilde{\Theta}_{12} \\ \tilde{\Theta}_{12}^\top & \tilde{\Theta}_{22} \ebm \bbm J & 0 \\ 0 & \hat{J}  \ebm x>0,
	\eeq
	for all $x \in \im (\Gamma) \backslash \{0\} \subset \R^{n+r+p}$ or, equivalently,
	$x \in \ker \bbm 0 &I_p & 0 & I_p \ebm \backslash \{0\}$. Let $R=\bbm 0 &I_p & 0 & I_p \ebm$, then by Finsler's lemma \cite{pinsler1936vorkommen}, \eqref{QMI_br_Gam} is equivalent to 
	\beq \nonumber 
	\bbm J & 0 \\ 0 & \hat{J} \ebm^\top  \bbm \tilde{\Theta}_{11} & \tilde{\Theta}_{12} \\ \tilde{\Theta}_{12}^\top & \tilde{\Theta}_{22} \ebm\bbm J & 0 \\ 0 & \hat{J} \ebm  -\mu R^\top R >0,
	\eeq
	for some $\mu$, which can be written as
	\beq \label{QMI_br}
	\bbm J  & 0 \\ 0 & \hat{J} \ebm^\top  \bbm \Theta_{11} & \Theta_{12} \\ \Theta_{12}^\top & \Theta_{22} \ebm\bbm J & 0 \\ 0 & \hat{J} \ebm>0,
	\eeq
	with 
	$$
	\Theta_{12}:= 
		\blkdiag( K_{12},- \mu I_p,-K_{12}, -\gamma^{-2}I_m)
	$$
	and 
	$$
	\Theta_{ii}:= 
	\blkdiag( K_{ii},
	(\frac{1}{2} -\mu)I_p,
    -K_{ii},
	-\gamma^{-2}I_m), 
	$$
	for $i=1,2$.
	The motivation of writing \eqref{e:br_in_QMI} in the form \eqref{QMI_br} is that the later form can be written in a QMI with the same quadratic variable as $J^\top N J\geq 0$. Namely, 
	by using a Schur complement argument, \eqref{QMI_br} is equivalent to
	\beq \label{e:QMI_J1}
	J^\top \left(\Theta_{11}-\Theta_{12}\hat{J} \left(\hat{J}^\top \Theta_{22} \hat{J}\right)^{-1}\hat{J}^\top \Theta_{12}^\top\right)J >0
	\eeq
	and $\hat{J}^\top \Theta_{22} \hat{J}>0$. This form allows us to use the matrix S-lemma in Proposition~\ref{matSlemma_strict} such that QMI ${J}^\top N {J} \geq 0$ implies \eqref{e:QMI_J1}. Particularly, \eqref{e:QMI_J1} holds with $J$ satisfying $J^\top {N} J  \geq0$ if and only if
	\beq \label{e:Slemma_T11}
	\Theta_{11}-\delta N-\Theta_{12}\hat{J} \left(\hat{J}^\top \Theta_{22} \hat{J}\right)^{-1}\hat{J}^\top \Theta_{12}^\top >0
	\eeq
	for some $\delta > 0$. To this end, we assume that $\hat{J}^\top \Theta_{22} \hat{J}>0$ holds for all $\hatJ$ satisfying $\hat{J}^\top N_{V,W} \hat{J} \geq 0$. We will see that this assumption is satisfied after completing the proof.
	
	Next, by using the (backward) Schur complement, \eqref{e:Slemma_T11} together with  $\hat{J}^\top \Theta_{22} \hat{J}>0$ is equivalent to 
	\beq \nonumber 
	\bbm 	\Theta_{11}-\delta N & \Theta_{12}\hat{J}\\
	\hat{J}^\top \Theta_{12}^\top &  \hat{J}^\top \Theta_{22} \hat{J} \ebm >0.
	\eeq 
	Then, a Schur complement with respect to the block matrix $\Theta_{11}-\delta N$ results in
	\beq \label{e;QMI_J2}
	\hat{J}^\top \left( \Theta_{22}-\Theta_{12}^\top \left(\Theta_{11}-\delta N\right)^{-1} \Theta_{12} \right) \hat{J} > 0 
	\eeq
	and $\Theta_{11}-\delta N >0$. Using Proposition~\ref{matSlemma_strict} again, \eqref{e;QMI_J2} holds for $\hat{J}$ satisfying $\hat{J}^\top {N}_{V,W} \hat{J}  \geq 0$ if and only if 
	$$
	\Theta_{22}-\eta N_{V,W}-\Theta_{12}^\top \left(\Theta_{11}-\delta N\right)^{-1} \Theta_{12} >0
	$$
	for some $\eta > 0$.  Finally, a (backward) Schur complement argument  implies that this is equivalent to \eqref{e:LMIs} as desired. Here, we have seen that $\hatJ$ satisfies $\hatJ^\top \Theta_{22} \hatJ >0$ for any $\hatJ$ satisfying $\hat{J}^\top {N}_{V,W} \hat{J}  \geq 0$ as an implication of \eqref{e;QMI_J2}.
\end{proof}


\subsection{Proof of Proposition~\ref{prop:upper_boundG0}} \label{app:boundsys0}
\begin{proof}
	Let ${\gamma}_0>0$. Then, from the bounded real lemma \cite{skelton1997unified,de1992discrete}, we have that ${\lVert \hat{\Sigma}_0 - \Sigma  \rVert_{\HHH_{\infty}} < {\gamma}_0}$ if and only if there exists  $K=K^\top>0$ in $\R^{(n+r)\times (n+r)}$ such that
	\beq  \label{e:br_lemma0}
\bbm K & 0 \\ 0 & I_p \ebm - 
\bbm A_{\mathrm{d}} & B_{\mathrm{d}} \\ C_{\mathrm{d}} & D_{\mathrm{d}} \ebm \bbm K & 0 \\ 0 & \gamma_0^{-2}I _m\ebm \bbm A_{\mathrm{d}} & B_{\mathrm{d}} \\ C_{\mathrm{d}} & D_{\mathrm{d}} \ebm^\top >0.
\eeq
	where
		\beq \nonumber
	A_\mathrm{d} \!:= \!\bbm A & 0 \\ 0 & \hatA_0 \ebm\!,  B_\mathrm{d} \!:= \!\bbm B \\ \hatB_0 \ebm, \ C_\mathrm{d} := \!\bbm C & -\hatC_0 \ebm, D_\mathrm{d} := \!D-\hatD_0.
	\eeq
	Next, we will show that condition \eqref{e:LMI_0} is equivalent to the existence of $K=K^\top>0$ such that \eqref{e:br_lemma0} holds for any system in $\sig$. 
	
First we introduce $J$ to denote the matrix as in \eqref{e:Jnotation}.
This allows to write \eqref{e:Ineq_noise} into $J^\top N J  \geq 0$ and moreover
\eqref{e:br_lemma0} into 
\beq \label{e:JTheta}
\bbm J^\top \bar{\Theta}_{11} J & J \bar{\Theta}_{12} \\ \bar{\Theta}_{12}^\top J & \bar{\Theta}_{22}
\ebm > 0,
\eeq
where 
\beq \nonumber 
 \bar{\Theta}_{11}\!:=\!\!\bbm 
K_{11} & 0  & 0  &  0    \cr
0& \! I_p\!-\!\hatC_0 K_{22} \hatC_0^\top \!-\! \gamma_0^{-2}\hatD_0\hatD_0^\top  & \hatC_0 K_{12}^\top & \gamma_0^{-2} \hatD_0  \cr
0 & K_{12}\hatC_0^\top & -K_{11} & 0  \cr
0 & \gamma_0^{-2}\hatD_0^\top  & 0  & -\gamma_0^{-2} I_m  \cr
\ebm\!,
\eeq 
\beq \nonumber 
\bar{\Theta}_{12}:= \bbm 
K_{12}   \cr
 \hatC_0 K_{22}\hatA_0^\top + \gamma_0^{-2} \hatD_0\hatB_0^\top  \cr
  -K_{12} \hatA_0^\top \cr
  -\gamma_0^{-2}\hatB_0^\top \cr
\ebm,
\eeq 
and $\bar{\Theta}_{22} := K_{22} - \hatA_0 K_{22} \hatA_0^\top - \gamma_0^{-2} \hatB_0 \hatB_0^\top$, which are denoting the block elements of the first matrix in \eqref{e:LMI_0}. Furthermore, the Schur complement of \eqref{e:JTheta} admits that $\bar{\Theta}_{22}>0$ and 
\beq \label{e:QMI_JthetaJ}
J^\top (\bar{\Theta}_{11}-\bar{\Theta}_{12}\bar{\Theta}_{22}^{-1}\bar{\Theta}_{12}^\top )J>0.
\eeq 
 Hence, by the strict matrix S-lemma in Proposition~\ref{matSlemma_strict}, \eqref{e:QMI_JthetaJ} holds with $J$ satisfying $J^\top N J \geq 0$ if and only if 
 \beq \label{e:QMI_JthetaSlemma}
 \bar{\Theta}_{11}-\bar{\Theta}_{12}\bar{\Theta}_{22}^{-1}\bar{\Theta}_{12}^\top - \delta N >0
 \eeq 
	for some $\delta >0$. Finally, \eqref{e:QMI_JthetaSlemma} and $\bar{\Theta}_{22}>0$ yield \eqref{e:LMI_0} via the (backward) Schur complement. 
\end{proof}
\bibliographystyle{IEEEtran}
\bibliography{IEEEabrv,reference.bib}
\end{document}